\documentclass{higherStructures}

\usepackage{float,tabularx}
\usepackage[utf8]{inputenc}
\usepackage{enumitem}   
\usepackage[english]{babel}
\usepackage{amsmath,amssymb,amsfonts,mathrsfs,mathtools,MnSymbol}

\usepackage{bbold,euscript}

\usepackage{quiver}

\usepackage{tikz}
\usetikzlibrary{arrows}
\usetikzlibrary{decorations.markings}

\title{Operadic Fibrations and Unary Operadic 2-categories}



\author{Dominik Trnka}

\address{Institute of Mathematics, University of Technology, Technick\'a 2896/2, 616 69, Brno, CZ.}

\email{trnka@fme.vutbr.cz}

\newtheorem{definition}[theorem]{Definition}
\newtheorem{proposition}[theorem]{Proposition}
\newtheorem{lemma}[theorem]{Lemma}
\newtheorem{corollary}[theorem]{Corollary}

\theoremstyle{definition}
\newtheorem{example}[theorem]{Example}
\newtheorem{remark}[theorem]{Remark}

\renewcommand{\epsilon}{\varepsilon}
\renewcommand{\phi}{\varphi}
\newcommand{\V}{\mathscr{V}}
\newcommand{\B}{\mathscr{B}}
\newcommand{\C}{\mathscr{C}}
\newcommand{\N}{\mathscr{N}}
\newcommand{\M}{\mathscr{M}}
\newcommand{\NC}{\mathscr{NC}}
\newcommand{\D}{\mathscr{D}}
\newcommand{\DC}{\mathscr{DC}}
\newcommand{\NDC}{\mathscr{NDC}}
\newcommand{\fib}{\mathrm{fib}}
\newcommand{\Cat}{\mathrm{Cat}}
\newcommand{\Set}{\mathrm{Set}}
\newcommand{\Hom}{\mathrm{Hom}}

\newcommand{\dec}{\mathrm{dec}_\top}
\newcommand{\StrMonCat}{\mathrm{StrMonCat}}
\newcommand{\Multi}{\mathrm{Multicat}}

\newcommand{\TwoUnOpCat}{\mathrm{uOpCat}}
\newcommand{\TwoCat}{\mathrm{Cat}}
\newcommand{\TwoCatNL}{\mathrm{Cat}_{\mathrm{NLax}}}

\newcommand{\sSet}{\mathrm{sSet}}

\newcommand{\sSetn}{\mathrm{sSet}_{\leq n}}
\newcommand{\sSett}{\mathrm{sSet}_{\leq 3}}

\newcommand\Bq[1]{\mathrm{Bq}^u(#1)}

\newcommand{\NN}{\mathbb{N}}
\newcommand{\uu}{\mathbb{1}}
\newcommand{\1}{\texttt{1}}

\newcommand{\bO}{\mathbb{O}}

\newcommand{\bP}{\mathbb{P}}
\newcommand{\bQ}{\mathbb{Q}}
\newcommand{\bM}{\mathbb{M}}
\newcommand{\tP}{\EuScript{P}}
\newcommand{\tM}{\EuScript{M}}

\newcommand{\sP}{\mathscr{P}}

\newcommand{\bC}{\mathbb{C}}

\newcommand{\ox}{\otimes}
\newcommand{\x}{\times}

\newcommand\card[1]{|#1|}

\newcommand\oper[2]{{#1}\textrm{-oper}(#2)}
\newcommand\coll[2]{{#1}\textrm{-coll}(#2)}

\def\colorop #1(#2;#3){
{#1}\left(
\rule{0pt}{15pt}\right.
	\hskip -3mm \begin{array}{c}
                 #3\\
		 #2
	 \end{array}
\hskip -3mm \left.
  \rule{0pt}{15pt} \right)
}

\def\tikzalpha (#1;#2){
\begin{scope}[shift={({#1},0)}]
\draw (0,0) --(0,1)node[circle,fill=black,inner sep=0pt,minimum size=6pt](p){};
 \draw  (-0.5,2) -- (p);
 \draw  (0.5,2) -- (p) ;
\draw (p) node[left]{$\alpha_{{#2}}$}; 
\end{scope}
}

\def\tikzf (#1){
\begin{scope}[shift={({#1},0)}]
\draw (0,0)--(0,2);
\draw (0,1);
\end{scope}
}

\def\tikzsigma (#1){
\begin{scope}[shift={({#1},0)}]
\draw (0,0)--(0,1)node[circle,fill=black,inner sep=0pt,minimum size=6pt](p){};
 \draw  (-1,2) -- (p);
 \draw  (0,2) -- (p) ;
 \draw  (1,2) -- (p) ;
\draw (p) node[right]{$\sigma_{0123}$}; 
\end{scope}
}

\def\comprule (#1;#2;#3;#4){
\begin{tikzpicture}
[thick, line cap=round,line join=round,
    x=1cm,y=1cm
    ] 
\draw (0,0)--(3,0)node[midway, above] {$#1$};
 \draw  (0,0) -- (1.5,-1)node[anchor=north east, midway] {$#2$};
 \draw  (1.5,-1) -- (3,0)node[midway, below right] {$#3$};
 \draw  (1.5,-0.3)node{$#4$};
\end{tikzpicture}
}

\DeclareMathOperator{\dom}{dom}
\DeclareMathOperator{\op}{op}
\DeclareMathOperator{\colim}{colim}

\let\pf\proof
\let\epf\endproof

\setlist[itemize]{leftmargin=*}

\begin{document}
\maketitle
\begin{abstract}
    We introduce unary operadic 2-categories as a framework for operadic Grothendieck construction of a categorical $\bO$-operad, $\bO$ being a unary operadic category. The construction is a fully faithful functor $\int_\bO$ which takes categorical $\bO$-operads to operadic functors over $\bO$, and we characterize its essential image by certain lifting properties. Such operadic functors are called operadic fibrations. Our theory is an extension of the discrete (unary) operadic case and, in some sense, of the classical Grothendieck construction of a categorical presheaf. For the terminal unary operadic category $\odot$, a categorical $\odot$-operad is a strict monoidal category $\V$ and its Grothendieck construction $\int_\odot\V$ is connected to the `para' construction appearing in machine learning. The 2-categorical setting provides a characterization of $\bO$-operads valued in $\V$ as operadic functors $\bO \to \int_\odot\V$. Last, we describe a left adjoint to $\int_\odot$. 
\end{abstract}
\begin{center}
\textit{In loving memory of grandfather Eduard (1940--2025).}
\end{center}

\section*{Introduction}
The theory of operadic categories was introduced in \cite{Batanin_Markl:duoidal} and further developed in \cite{Batanin_Markl:kodu2022, Batanin_Markl:kodu2021}, and \cite{Batanin_Markl:blob}, as a unifying framework for compositional structures and their homotopy theory. The last of the above mentioned works focuses only on `unary' operadic categories which themselves describe a lot of known structures. On one hand, every discrete decomposition space is a unary operadic category \cite{GKW}. On the other, operads of certain unary operadic categories are monoids, categories, sharades, \&c. 
In~\cite{Batanin_Markl:blob}, unary operadic categories were used for a characterization of the `blob complex' of \cite{MW} as a bar construction of a certain operad. The~unary case often serves as a toy model for new directions in the theory of (non-unary) operadic categories, cf.~\cite{bivar,Hack}.

It is well known that the Grothendieck construction establishes a correspondence of categorical presheaves and categorical fibrations.
In the operadic context, the Grothendieck construction for $\Set$-valued operads, as well as discrete operadic fibrations, were introduced already in the first paper \cite{Batanin_Markl:duoidal}. The construction is a part of an equivalence 
$$\oper{\bO}{\Set} \simeq \textrm{DoFib}(\bO)$$
 between $\Set$-valued $\bO$-operads discrete operadic fibrations over an operadic category $\bO$. 
Our goal is to establish a correspondence 
\begin{equation}\label{equation:correspondence}\oper{\bO}{\Cat} \simeq \textrm{oFib}(\bO)
\end{equation}
 between categorical $\bO$-operads and operadic functors into $\bO$ with certain lifting properties, which we call operadic fibrations. The lifting properties are operadic analogues of the properties of classical categorical fibrations.
 
\begin{figure}[H]
    \centering
   \begin{tikzpicture}
[very thick, line cap=round,line join=round,
    x=0.6cm,y=0.4cm
    ] 
\draw  (2.,8.)  node[circle,fill=black,inner sep=0pt,minimum size=5pt] {} -- (1.,6.);
\draw  (2.,8.)-- (2.,6.);
\draw  (2.,8.)-- (3.,6.);
\draw  (1.,6.)-- (1.,4.);
\draw  (2.,6.)-- (2.,4.);
\draw  (3.,6.)-- (3.,4.);
\draw  (3.,6.)-- (4.,4.);
\draw  (1.,4.)-- (0.,2.);
\draw  (2.,4.)-- (1.,2.);
\draw  (2.,4.)-- (2.,2.);
\draw  (2.,4.)-- (3.,2.);
\draw  (3.,4.)-- (4.,2.);
\draw  (4.,4.)-- (5.,2.);
\draw  (0.,2.)-- (0.,0.)  node[circle,fill=black,inner sep=0pt,minimum size=5pt] {};
\draw  (1.,2.)-- (1.,0.) node[circle,fill=black,inner sep=0pt,minimum size=5pt] {};
\draw  (2.,2.)-- (2.,0.) node[circle,fill=black,inner sep=0pt,minimum size=5pt] {};
\draw  (3.,2.)-- (3.,0.) node[circle,fill=black,inner sep=0pt,minimum size=5pt] {};
\draw  (4.,2.)-- (4.,0.) node[circle,fill=black,inner sep=0pt,minimum size=5pt] {};
\draw  (4.,2.)-- (5.,0.) node[circle,fill=black,inner sep=0pt,minimum size=5pt] {};
\draw  (4.,2.)-- (6.,0.) node[circle,fill=black,inner sep=0pt,minimum size=5pt] {};
\draw  (5.,2.)-- (7.,0.) node[circle,fill=black,inner sep=0pt,minimum size=5pt] {};
\draw [thick, gray] (-1.,2.)-- (7.,2.);
\draw [thick, gray] (-0.62,4.)-- (5.56,4.);
\draw[thick, gray]  (-0.82,6.)-- (6.,6.);
\end{tikzpicture}
\begin{tikzpicture}
[very thick, line cap=round,line join=round,
    x=0.6cm,y=0.4cm
    ] 
\draw  (2.,8.) node[circle,fill=black,inner sep=0pt,minimum size=5pt] {}-- (1.,6.);
\draw  (2.,8.)-- (2.,6.);
\draw  (2.,8.)-- (3.,6.);
\draw  (4.,4.)-- (5.,2.);
\draw  (0.,2.)-- (0.,0.);
\draw  (1.,2.)-- (1.,0.);
\draw  (2.,2.)-- (2.,0.);
\draw  (3.,2.)-- (3.,0.);
\draw  (4.,2.)-- (4.,0.);
\draw  (4.,2.)-- (5.,0.);
\draw  (4.,2.)-- (6.,0.);
\draw  (5.,2.)-- (7.,0.);
\draw  (4.,4.)-- (4.,2.);
\draw  (3.,4.)-- (3.,2.);
\draw  (2.,4.)-- (2.,2.);
\draw  (1.,4.)-- (1.,2.);
\draw  (0.,2.)-- (0.,4.);
\draw  (0.,4.)-- (1.,6.);
\draw  (2.,6.)-- (1.,4.);
\draw  (2.,4.)-- (2.,6.);
\draw  (2.,6.)-- (3.,4.);
\draw  (3.,6.)-- (4.,4.);

\draw  (0.,2.)-- (0.,0.)  node[circle,fill=black,inner sep=0pt,minimum size=5pt] {};
\draw  (1.,2.)-- (1.,0.) node[circle,fill=black,inner sep=0pt,minimum size=5pt] {};
\draw  (2.,2.)-- (2.,0.) node[circle,fill=black,inner sep=0pt,minimum size=5pt] {};
\draw  (3.,2.)-- (3.,0.) node[circle,fill=black,inner sep=0pt,minimum size=5pt] {};
\draw  (4.,2.)-- (4.,0.) node[circle,fill=black,inner sep=0pt,minimum size=5pt] {};
\draw  (4.,2.)-- (5.,0.) node[circle,fill=black,inner sep=0pt,minimum size=5pt] {};
\draw  (4.,2.)-- (6.,0.) node[circle,fill=black,inner sep=0pt,minimum size=5pt] {};
\draw  (5.,2.)-- (7.,0.) node[circle,fill=black,inner sep=0pt,minimum size=5pt] {};
\draw [thick, gray] (-1.,2.)-- (7.,2.);
\draw [thick, gray] (-0.62,4.)-- (5.56,4.);
\draw[thick, gray]  (-0.82,6.)-- (6.,6.);
\end{tikzpicture}
\begin{tikzpicture}
[very thick, line cap=round,line join=round,
    x=0.6cm,y=0.4cm
    ] 

\draw  (2.,8.) node[circle,fill=black,inner sep=0pt,minimum size=5pt] {}-- (1.,6.);
\draw  (2.,8.)-- (2.,6.);
\draw  (2.,8.)-- (3.,6.);
\draw  (3.,6.)-- (3.,4.);
\draw  (3.,6.)-- (4.,4.);
\draw  (3.,4.)-- (4.,2.);
\draw  (3.,4.)-- (3.,2.);
\draw  (2.,6.)-- (2.,4.);
\draw  (1.,6.)-- (1.,4.);
\draw  (3.,2.)-- (4.,0.) node[circle,fill=black,inner sep=0pt,minimum size=5pt] {};
\draw  (4.,2.)-- (5.,0.) node[circle,fill=black,inner sep=0pt,minimum size=5pt] {};
\draw  (4.,4.)-- (6.,2.);
\draw  (3.,4.)-- (5.,2.);
\draw  (5.,2.)-- (6.,0.) node[circle,fill=black,inner sep=0pt,minimum size=5pt] {};
\draw  (6.,2.)-- (7.,0.) node[circle,fill=black,inner sep=0pt,minimum size=5pt] {};
\draw  (2.,2.)-- (3.,0.) node[circle,fill=black,inner sep=0pt,minimum size=5pt] {};
\draw  (2.,2.)-- (2.,0.) node[circle,fill=black,inner sep=0pt,minimum size=5pt] {};
\draw  (2.,2.)-- (1.,0.) node[circle,fill=black,inner sep=0pt,minimum size=5pt] {};
\draw  (2.,4.)-- (2.,2.);
\draw  (1.,4.)-- (1.,2.);
\draw  (1.,2.)-- (0.,0.) node[circle,fill=black,inner sep=0pt,minimum size=5pt] {};

\draw [thick, gray] (-1.,2.)-- (7.,2.);
\draw [thick, gray] (-0.62,4.)-- (5.56,4.);
\draw[thick, gray]  (-0.82,6.)-- (6.,6.);
\end{tikzpicture}
    \caption{Three isomorphic leveled trees.}
    \label{fig:leveled_trees}
\end{figure}
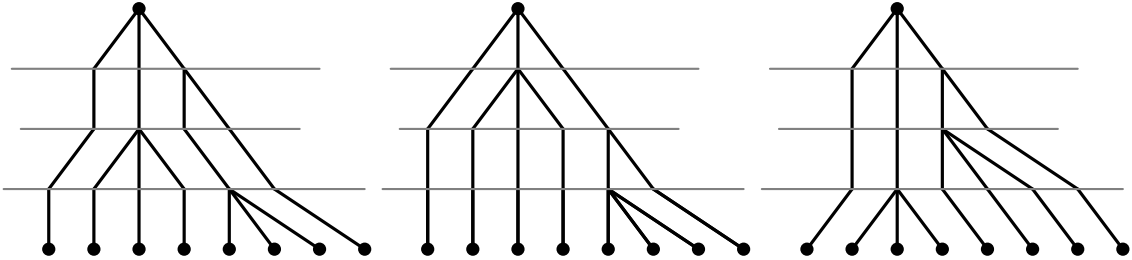

 Categorical operads comprise operads with values in groupoids and partially ordered sets. A~remarkable application of operads valued in groupoids is \cite{Fresse} on Grothendieck--Teichmüller groups. Operads valued in posets appear e.g.~in \cite{Berger97} and many poset-valued operads are described in \cite{Bashkirov}. Categorical operads are further considered, for instance, in \cite{Batanin08,CG14}.
A concrete example of a (non-strict) categorical operad is the operad of leveled trees. The operations in arity $n$ are leveled planar rooted trees with $n$ leaves, and there is a unique isomorphism of trees if they differ by an admissible change of leveling, cf.~Figure~\ref{fig:leveled_trees}. This operad plays a key r{\^o}le in the construction of free operads~\cite{Batanin_Markl:kodu2021,moje}. 

 It turns out that the framework of operadic \hbox{(1-)categories} is not sufficient for the correspondence \eqref{equation:correspondence}, and that a 2-categorical structure is needed. The 2-categorical setting brings also other benefits, namely, an identification of $\bO$-operads with operadic functors out of~$\bO$. This interprets $\bO$-operads as generalized presheaves on $\bO$. As a result, $\bO$-operads are defined intrinsically to the theory of operadic categories. Last, a categorical operad of the terminal operadic category is a strict monoidal category and its operadic Grothendieck construction is related to the `para' construction which appears in machine learning, for example in~\cite{para1}.
 
This paper derives the results in the unary case only. An extension of the theory to general operadic categories can be expected in our future work.
 Before describing the results in more detail, we give an intuitive idea about unary operadic categories, their operads, and the discrete operadic Grothendieck construction.
 
\subsection*{Shapes and structures.} This is an informal introduction to
unary operadic categories which are recalled in Definition~\ref{definition:unary_operadic_category}. An operadic category is a structure suitable for handling shapes and composionality. Hence, in~this motivational subsection we shall use a more suggestive name \textit{shape calculus} for unary operadic categories. A shape calculus $\bO$ is specified by a class of \textit{shapes} (objects of $\bO$) and triangular \textit{composition rules} (morphisms of $\bO$) of the form 

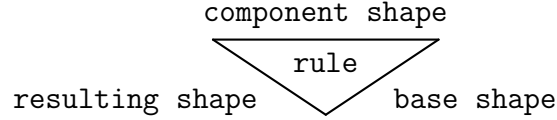
\[\comprule (\texttt{component shape};\texttt{resulting shape};\texttt{base shape};\texttt{rule})\]
A rule consists of three shapes and an information of how to compose the shapes. Additionally to source and target, every rule
$$\texttt{resulting shape}\xrightarrow{\texttt{rule}}\texttt{base shape}$$
of $\bO$ also has a \textit{fiber}, the \texttt{component shape}. The fiber assignment is assembled into functors $\phi_s\colon\bO/s\to\bO$ from the slice category, for all shapes $s$. Lastly, there are \textit{trivial shapes} and two trivial composition rules 
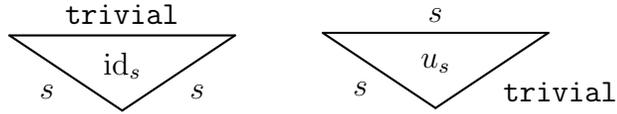
\[\hspace{1.5cm}
\begin{tikzpicture}
[thick, line cap=round,line join=round,
    x=1cm,y=1cm
    ] 
\draw (0,0)--(3,0)node[midway, above] {$\texttt{trivial}$};
 \draw  (0,0) -- (1.5,-1)node[anchor=north east, midway] {$s$};
 \draw  (1.5,-1) -- (3,0)node[midway, below right] {$s$};
 \draw  (1.5,-0.4)node{$\mathrm{id}_s$};
\end{tikzpicture}
\hspace{1cm}\
\begin{tikzpicture}
[thick, line cap=round,line join=round,
    x=1cm,y=1cm
    ] 
\draw (0,0)--(3,0)node[midway, above] {$s$};
 \draw  (0,0) -- (1.5,-1)node[anchor=north east, midway] {$s$};
 \draw  (1.5,-1) -- (3,0)node[midway, below right] {$\texttt{trivial}$};
 \draw  (1.5,-0.4)node{$u_s$};
\end{tikzpicture}\]
 for every shape $s$, where both have a meaning of `do nothing with $s$'. These correspond to an identity map and to a special map into a trivial shape in $\bO$.

 For a calculus $\bO$, an \textit{$\bO$-operad}  $\tP$ (cf.~Definition~\ref{definition:O-operad_classical}) is a compositional structure which consists of sets $\tP_s$ for every shape $s \in \bO$, and composition maps $$\mu_\texttt{rule}\colon \tP_\texttt{base} \x \tP_\texttt{comp} \to \tP_\texttt{result}$$ indexed by the composition rules. These maps are assumed to be unital and associative, which is expressed in terms of trivial objects and axioms of operadic categories, which we do not state here. 
 Four examples are as follows.

\begin{example}\label{example:trivial_calculus}
     The trivial calculus $\odot$ consists only of one trivial shape $t$ and one trivial composition rule. \[
\begin{tikzpicture}
[thick, line cap=round,line join=round,
    x=1cm,y=1cm
    ] 
\draw (0,0)--(3,0)node[midway, above] {$t$};
 \draw  (0,0) -- (1.5,-1)node[anchor=north east, midway] {$t$};
 \draw  (1.5,-1) -- (3,0)node[midway, below right] {$t$};
 \draw  (1.5,-0.4)node{$\mathrm{id}_t$};
\end{tikzpicture}\]
An $\odot$-operad is a single set $\tP_{t}=M$, together with a single map $M\x M\xrightarrow{\cdot} M$, which is associative and unital. The $\odot$-structures are precisely monoids. 
\end{example}
\begin{example}\label{example:calculus_of_naturals}
     Let $(\mathbb{N},+,0)$ be the monoid of natural numbers. The shapes of calculus~$\NN_\triangle$ are natural numbers and the composition rules are 
     \[\begin{tikzpicture}
[thick, line cap=round,line join=round,
    x=1cm,y=1cm
    ] 
\draw (0,0)--(3,0)node[midway, above] {$k$};
 \draw  (0,0) -- (1.5,-1)node[anchor=north east, midway] {$n$};
 \draw  (1.5,-1) -- (3,0)node[midway, below right] {$m$};
 \draw  (1.5,-0.2)node{$m+k=n$};
\end{tikzpicture}
\]
whenever $m+k=n$. An $\NN_\triangle$-operad is a collection of sets $\tP_n$ and composition maps 
$$\tP_m\x\tP_k\xrightarrow{\mu_{m,k}}\tP_{m+k},$$
which are associative and unital. This is precisely a unital graded algebra (a graded monoid). More generally, any monoid $M$ determines a calculus $M_{\triangle}$ and $M_{\triangle}$-operads are
lax monoidal functors $M\to (\Set,\times,\1)$, where $M$ is regarded as a discrete monoidal category.
\end{example}
\begin{example}\label{example:calculus_bouquets}
    For a fixed set $S$, $\Bq{S}$ stands for a calculus of `unary bouquets'. Its~shapes are pairs $(a,b)$ of elements of $S$, written as $a\to b$. The composition rules are \[\comprule (a\to b;a\to c';b\to c'';{c'=c''})\] whenever $c'=c''$. A $\Bq{S}$-operad is a collection of sets $\tP_{a\to b}=\Hom(a,b)$ and compositions $$\Hom(b,c)\x \Hom(a,b) \to \Hom(a,c).$$
    The $\Bq{S}$-operads are small categories with set of objects $S$. The name \textit{bouquets} is more meaningful in the non-unary case, cf.~\cite{Batanin_Markl:duoidal}.
    \end{example}
\begin{example}\label{example:calculus_subsets}
For a fixed set $X$, the shapes of calculus $\mathrm{sub}(X)$ are subsets of $X$. The~composition rules are \[\comprule (B-A;B;A;{A\subseteq B})\] whenever ${A\subseteq B}$. The $\mathrm{sub}(X)$-structures are structures describing, for example, gluing of local functions on the set $X$. Typically the set $X$ has some additional structure (algebraic, topological, combinatorial, \&c.), cf.~\cite[Part 2]{Batanin_Markl:blob}.
\end{example}

\subsection*{Discrete operadic Grothendieck construction.}
There is an interplay between the shape calculi and their operads. Having an $\bO$-operad $\tP$, its operadic Grothendieck construction $\int \tP$ produces a new calculus, whose shapes are all elements of the operad $\tP$ and the composition rules are given by the composition maps $\mu$ of $\tP$. We have already seen an example of this phenomenon. Namely, a monoid $M$ is an $\odot$-operad and the construction produces the calculus $\int M=M_\triangle$. 

Formally, the shapes of calculus $\int\tP$ are pairs $(s,\alpha)$, where $s$ is a shape of $\bO$ and $\alpha \in \tP_s$. The rules of $\int \tP$ have form 
\[
\begin{tikzpicture}
[thick, line cap=round,line join=round,
    x=1cm,y=1cm
    ] 
\draw (0,0)--(6,0)node[midway, above] {$(s,\alpha)$};
 \draw  (0,0) -- (3,-2)node[anchor=north east, midway] {$(p,\gamma)$};
 \draw  (3,-2) -- (6,0)node[midway, below right] {$(q,\beta)$};
 \draw  (3,-0.6)node{$\mu_{r}(\beta,\alpha)=\gamma$};
\end{tikzpicture}
\]
whenever $\mu_{r}(\beta,\alpha)=\gamma$, where $r$ is a rule of $\bO$ with 
\[\begin{tikzpicture}
[thick, line cap=round,line join=round,
    x=1cm,y=1cm
    ] 
\draw (0,0)--(3,0)node[midway, above] {$s$};
 \draw  (0,0) -- (1.5,-1)node[anchor=north east, midway] {$p$};
 \draw  (1.5,-1) -- (3,0)node[midway, below right] {$q$};
 \draw  (1.5,-0.4)node{$r$};
\end{tikzpicture} \]
and $\alpha \in \tP_s$, $\beta \in \tP_q$ and $\gamma \in \tP_p$. 
The associativity of $\tP$ translates into associativity of composition of the rules of $\int\tP$, and the units of $\tP$ become trivial objects. The (two-sided) unitality of~$\tP$ provides identities for composition in $\int \tP$ as well as the second trivial rule. Further, there is a direct comparison (an operadic functor) of the calculi $\int\tP $ and $ \bO$ which is the projection $\pi\colon \int\tP \to \bO$ on the first component. 

The projection $\pi$ has a special property: given any two shapes $S$ and $Q$ of $\int\tP$, and any rule~$r$ in $\bO$ with base shape $\pi(Q)=q$, component $\pi(S)=s$ and resulting shape $p$, there is a unique rule $\Tilde{r}$ having component $S$, base $Q$ and result $P$, such that $\pi(P)=p$. The property is drawn in the following diagram. 
\[
\begin{tikzpicture}
[thick, line cap=round,line join=round,
    x=1cm,y=1cm
    ] 
\draw (0,0)--(6,0)node[midway, above] {$S$};
 \draw[dotted]   (0,0) -- (3,-2)node[midway,below left,draw,circle,dotted] {$P$};
 \draw  (3,-2) -- (6,0)node[midway, below right] {$Q$};
 \draw  (3,-0.8)node[draw,circle,dotted](y){$\Tilde{r}$};
\begin{scope}[shift={(1.5,-3)}]
    \draw (0,0)--(3,0)node[midway, above] {$s$};
 \draw  (0,0) -- (1.5,-1)node[anchor=north east, midway] {$p$};
 \draw  (1.5,-1) -- (3,0)node[midway, below right] {$q$};
 \draw  (1.5,-0.4)node[circle] (x){$r$};
 \draw[|->,out=130,in=-130] (x) to node[midway,left] {$\exists!$} (y);
\end{scope}
\end{tikzpicture} \] 
 We say that $\pi$ is a \textit{discrete operadic fibration}. 
The unique lifting property characterizes the calculi of the form $\int\tP$ and it allows us to compare $\bO$-operads and $\int\tP$-operads.
A natural question is what is the construction for operads valued in categories, and what are its properties. We now describe our results in more detail.

\subsection*{Main results.}
 We introduce Definition~\ref{definition:unary_operadic_2-category_A} of a unary operadic 2-category, which is motivated by the characterization of unary operadic categories of \cite{GKW} via simplicial sets and décalage, recalled in Proposition~\ref{proposition:dec_and_nerve}. We reformulate the definition in terms of (lax normalized) fiber functors and a choice of local (lali-)terminal objects, which resembles the original definition of \cite{Batanin_Markl:duoidal}, cf.~Definition~\ref{definition:unary_operadic_category}.
 In Definition \ref{definition:operad2} we introduce operads of unary operadic 2-categories and relate them to operadic functors $\tP\colon \bO \to \V_\triangle$, where the target is the operadic 2-category of lax triangles in a monoidal category $\V$. For $\V=\Cat$, we imitate techniques known for presheaves and define the operadic Grothendieck construction of a (1-connected) categorical operad as a certain pullback, cf.~Definition~\ref{definition:G-construction}. The extension to non-1-connected operads is immediate and the construction is a fully faithful functor 
 $$\int_\bO\colon \oper{\bO}{\Cat}\to\TwoUnOpCat/\bO$$
into the category of unary operadic 2-categories over $\bO$. Our main result is Theorem~\ref{theorem:fibrations}, which characterizes the image of $\int_\bO$ as operadic functors into $\bO$ with certain lifting properties.
The~properties are operadic analogues of the lifting properties of classical categorical split fibrations and we call the functors split operadic fibrations, cf.~Definition~\ref{definition:fibration}. Lastly, in the special case of $\bO=\odot$, the terminal unary operadic category, Section~\ref{section:left_adjoint_*} describes a left adjoint to $\int_\bO$, which gives a comparison of unary operadic 2-categories to strict monoidal categories:
\[\begin{tikzcd}
            \StrMonCat \arrow[r, shift left=1ex, "\int"{name=U}] & \TwoUnOpCat. \arrow[l, shift left=.5ex, "\Phi"{name=F}]
            \arrow[phantom, from=F, to=U, , "\scriptscriptstyle\boldsymbol{\top}"]
        \end{tikzcd}\]
The adjunction follows from an adjunction 
\[\begin{tikzcd}
            \StrMonCat \arrow[r, shift left=1ex, ""{name=U}] & \sSet_{\leq3}, \arrow[l, shift left=.5ex, ""{name=F}]
            \arrow[phantom, from=F, to=U, , "\scriptscriptstyle\boldsymbol{\top}"]
        \end{tikzcd}\]
which might be of an independent interest.

\section{Unary operadic categories and their operads}
This section recalls the definition of a (unital) unary operadic category~$\bO$, the definition of a (unital) $\bO$-operad, operadic Grothendieck construction for set-values operads, and discrete operadic fibrations. The material is mostly taken from \cite{Batanin_Markl:blob}. We also recall the result of \cite{GKW}, which characterises unary operadic categories using simplicial sets and décalage construction, and we include some examples.

Let us denote by $\uu_x$ the identity morphism on an object~$x$. For a category $\C$ and its object $x$, let $\C/x$ denote the slice category over $x$. Observe that any functor $F\colon \C \to \D$ induces a functor on slices $F/x\colon \C/x \to \D/F(x)$. We denote by $\Set$ the category of sets and $\Cat$ stands for the category of (small) categories.

\begin{definition}\label{definition:unary_operadic_category}
A unary operadic category is a category $\bO$ together with a family of functors $ \fib_x\colon \bO/x \to \bO$, indexed by objects $x$ of $\bO$, such that for any map $f\colon y \to x$ of $\bO$, the following diagram commutes.
\begin{equation}\label{diagram:fiber_axiom}
\begin{tikzcd}
	{(\bO/x)/f} & {\bO/\fib_x(f)} \\
	{\bO/y} & {\bO}.
	\arrow["{\dom/f}"', from=1-1, to=2-1]
	\arrow["{\fib_y}"', from=2-1, to=2-2]
	\arrow["{\fib_x/f}", from=1-1, to=1-2]
	\arrow["{\fib_{\fib_x(f)}}", from=1-2, to=2-2]
\end{tikzcd}\end{equation}
We further require for every connected component $c$ of $\bO$ a choice of local terminal object~$u_c$, which belongs to the component $c$, such that $\fib_{u_c}\colon \bO/u_c \to \bO$ is the domain functor, and for every object $x$ of $\bO$, $\fib_x(\uu_x)$ is one of the objects $u_c$. 
\end{definition}

The objects $u_c$ are called \textit{trivial} and the functors $\fib_x$ are called \textit{fiber} functors. For an object $f\colon y\to x$ of $\bO/x$, the value $\fib_x(f)$ is called the \textit{fiber} of $f$ and it will be denoted simply by $\fib(f)$. For a morphism $g\colon f\circ g\to f$ in $\bO/x$, the induced map between fibers $\fib_x(g)\colon \fib(f\circ g)\to \fib(f)$ will be denoted by $g_{f}$. A \textit{quasibijection} is a morphism whose fiber is trivial. The fact that an object $r$ is a fiber of a morphism $f\colon y\to x$ is usually expressed by $$r\triangleright y \xrightarrow{f} x.$$ From the simplicial viewpoint on unary operadic categories, recalled in Proposition~\ref{proposition:GKW} below, one can view the map $f$ as a triangle 
\[\begin{tikzcd}
	\cdot && \cdot \\
	& \cdot
	\arrow[""{name=0, anchor=center, inner sep=0},"{r}", from=1-1, to=1-3]
	\arrow["{y}"', from=1-1, to=2-2]
	\arrow["{x}", from=1-3, to=2-2]
	\arrow["f"{description},draw=none, from=0, to=2-2]
	\end{tikzcd}\]
    This explains the use of triangular rules in the introduction.
\begin{definition}
An operadic functor between unary operadic categories is a functor which respects the fiber functors and preserves trivial objects.
\end{definition}

\begin{definition}\label{definition:O-operad_classical}
Let $\bO$ be a unary operadic category and $\V$ a monoidal category. A~(unital) $\bO$-operad $\tP$ with values in $\V$ is a collection of objects $\tP_x \in \V$, indexed by objects of~$\bO$, together with composition maps
$$\tP_x \ox \tP_r \xrightarrow{\tP_f} \tP_y,$$
for every $r \triangleright y \xrightarrow{f} x$ in $\bO$, which are associative and unital in the following sense. For~any composable maps $r \triangleright y \xrightarrow{f} x$ and $s \triangleright z \xrightarrow{g} y$, with $t \triangleright z \xrightarrow{f\circ g} x$ and the induced map $s \triangleright r \xrightarrow{g_f} t$, the diagram
\[\begin{tikzcd}
	{\tP_x \ox \tP_r \ox \tP_s} & {\tP_y\ox \tP_s} \\
	{\tP_x \ox \tP_t} & {\tP_z}
	\arrow["{\tP_f \ox \tP_s }", from=1-1, to=1-2]
	\arrow["{\tP_{g}}", from=1-2, to=2-2]
	\arrow["{\tP_{f\circ g}}"', from=2-1, to=2-2]
    \arrow["{\tP_x\ox \tP_{g_f}}"', from=1-1, to=2-1]
\end{tikzcd}
\]
commutes.
For every connected component $c$ of $\bO$, the operad $\tP$ is equipped with unit maps
$$\epsilon_{c}\colon \1 \to \tP_{u_c},$$
from the monoidal unit $\1$ of $\V$, such that 
for any object $x$ of $\bO$ with $\fib(\uu_x)=u_c$, the diagram
\[\begin{tikzcd}
	{\tP_x} & {\tP_x} \\
	{ \tP_x \ox \1} & {\tP_x\ox \tP_{u_c}}
	\arrow["\cong"', from=1-1, to=2-1]
	\arrow["{ \tP_x\ox \epsilon_c}"', from=2-1, to=2-2]
	\arrow["{\tP_{\uu_x}}"', from=2-2, to=1-2]
	\arrow["\uu", from=1-1, to=1-2]
\end{tikzcd}\]
commutes, and for any $x\xrightarrow{!_d} u_d$, the unique map to one of the chosen local terminals, the diagram
\[\begin{tikzcd}
	{\tP_x} & {\tP_x} \\
	{\1\ox \tP_x} & {\tP_{u_d}\ox \tP_x}
	\arrow["\cong"', from=1-1, to=2-1]
	\arrow["{\epsilon_d\ox \tP_x}"', from=2-1, to=2-2]
	\arrow["{\tP_{!_d}}"', from=2-2, to=1-2]
	\arrow["\uu", from=1-1, to=1-2]
\end{tikzcd}\]
commutes. 
\end{definition}

When the objects of the category $\V$ are concrete (e.g.~when $\V$ is $\Set$ or $\Cat$), for elements $a\in\tP_x$, $b\in\tP_r$ we write $a\cdot_{f}b$ instead of $\tP_f(a,b)$. The above diagrams translate as  
$$(a\cdot_{f} b)\cdot_{g} c = a\cdot_{f\circ g}(b\cdot_{g_f} c),$$
 $$ a\cdot_{\uu_x}e_c =a,$$
  $$ e_c \cdot_! b= b,$$
where $e_c$ is the image of $\epsilon_c$. 
\begin{example}\label{example:odot}
Let $\odot$ be the terminal category with one object $*$ and one morphism~$\uu_{*}$. It~is an operadic category with $* \triangleright * \xrightarrow{\uu_{*}} *$.
Analogously to~\cite[Exercise~39]{Batanin_Markl:blob}, $\odot$-operads in $\Cat$ are strict monoidal categories.
\end{example}

\begin{example}
Let $S$ be a set. The operadic category $\Bq{S}$ of unary bouquets has as objects pairs $(a,b)$ of elements of $A$. There is a unique morphism $(a',b')\to (a'',b'')$ whenever $b'=b''$, with fiber $(a',a'')$. Consequently, $\Bq{S}$-operads with values in $\Set$ are small categories whose set of objects is $S$, and $\Bq{S}$-operads valued in $\Cat$ are (strict) 2-categories with the set of objects~$S$. 
\end{example}

\begin{example}\label{example:O_M}
Let $(\M,\cdot,e)$ be a strict monoidal category. We construct the category~$\M_\triangle$ as follows. The objects of $\M_\triangle$ are the same as objects of $\M$, but a morphisms $y\to x$ of $\M_\triangle$ is a pair $(a,\alpha)$, where $a \in M$ and $\alpha\colon x \cdot a \to y$ is a map in~$\M$. The identity on an object $x$ is the pair $(e,\uu_x)$ and the composite of two morphisms
$$ z \xrightarrow{(a, \alpha)} y \xrightarrow{(b, \beta)} x$$
is $\big(b \cdot a, \alpha\circ (\beta\cdot a)\big)$:
$$
 x\cdot b\cdot a \xrightarrow{\beta\cdot a} y\cdot a \xrightarrow{\alpha}z.
$$
If we define the fiber of a map $(a,\alpha)$ to be the object $a$ and for 
$$ z \xrightarrow{(a, \alpha)} y \xrightarrow{(b, \beta)} x,$$
the induced fiber map $(a,\alpha)_{(b,\beta)}$ to be $(a,\uu_{b\cdot a})$, we \textit{almost} get a unary operadic category. The problem is, that there are in general no local terminal objects. Moreover, the fiber functor forgets the information of the morphism $\alpha$. We shall resolve the problems by introducing unary operadic 2-categories. When passing to the 2-categorical framework, instead of being terminal, the monoidal unit $e \in \M_\triangle$ has the property, that for any object $x$ of $\M$, the hom-category $\M_\triangle(x,e)$ has a terminal object, namely the map $(x,\uu_x)$. More generally, the maps $(a,\uu_{b\cdot a})$ will play an important r{\^o}le in the theory. In machine learning the category $\M_\triangle$ is known as the `para' construction of a monoidal category \cite{para1,para2,para3}.
\end{example}
\begin{definition}
Let $p\colon \bP \to \bO$ be an operadic functor between two unary operadic categories. It is a discrete operadic fibration, if $p$ induces an epimorphism $\pi_0(\bP) \to \pi_0(\bO)$ on connected components, and for any morphism $f\colon y\to x$ in $\bO$ with fiber $r$, and any two objects $R,X$ of~$\bP$, such that
$p(X) = x$ and $p(R) = r$, there exists a unique object $Y$ and a unique morphism $F \colon Y \to X$ in $\bP$ with fiber $R$,
such that $p(F) = f$.
\end{definition}

\begin{definition}
Let $\bO$ be a unary operadic category and $\tP$ an $\bO$-operad in $\Set$. Its~operadic Grothendieck construction is a unary operadic category $\int_\bO \tP$ whose objects are pairs $(x,a)$, where~$x$ is an object of $\bO$ and $a\in \tP_x$. A morphisms $(x',a')\to (x'',a'')$ of $\int_\bO \tP$ is a pair $(f,c)$ of a map $f\colon x' \to x''$ and element $c\in \tP_{\fib(f)}$, such that $a'=a''\cdot_f c$.    
\end{definition}
The operadic structure of $\int_\bO \tP$ is described in more detail in \cite{Batanin_Markl:blob}. We~recall, that the projection $\pi_\tP\colon \int_\bO \tP \to \bO$, sending a pair $(x,a)$ to $x$, is a discrete operadic fibration. 
\begin{proposition}{\cite[Proposition 35]{Batanin_Markl:blob}}
The~assignment $\tP\mapsto \pi_\tP$ is part of an equivalence
$$\oper{\bO}{Set}\simeq \mathrm{DoFib}(\bO)$$
between the category of set-valued $\bO$-operads and discrete operadic fibrations over $\bO$.
\end{proposition}
 \begin{example}
For a discrete monoidal category $\M$ (equivalently $\M$ is a monoid in $\Set$), the only maps of the category $\M_\triangle$ of Example~\ref{example:O_M} are $b\cdot a \xrightarrow{(a,\uu_b\cdot a)} b$, and there is an isomorphism $\M_\triangle\cong \int_\odot \M$, where $\M$ is regarded as an $\odot$-operad in $\Set$. Hence the unique map $\M_\triangle\to \odot$ is a discrete operadic fibration.
\end{example}
We close this section by recalling a characterisation of unary operadic categories via simplicial sets and décalage construction.
\begin{definition}\label{definition:dec}
    An upper décalage $\dec X$ of simplicial set $X$ is a simplicial set
    $$(\dec X)_n = X_{n+1},$$
    without the top face and degeneracy maps of $X$, i.e.~$d_i^{\mathrm{dec}}=d_i$ and $s_j^{\mathrm{dec}}=s_j$.  
\end{definition}
For every category $\C$, its nerve $\N\C$ is a simplicial set of sequences of composable morphisms. The assignment $\N$ is a fully faithful functor $\Cat \to \sSet$ and the simplicial set $\N\C$ is 2-coskeletal. This means that any element of $\NC_{k>2}$ can be expressed as a tuple of compatible elements $t_1,\cdots,t_k \in \NC_2$, and thus the simplicial set $\NC$ is determined by the data
\[\begin{tikzcd}
	{\NC_0} & {\NC_1} & {\NC_2.}
	\arrow["{s_0}"{description}, from=1-1, to=1-2]
	\arrow["{d_1}"', shift right=3, from=1-2, to=1-1]
	\arrow["{d_0}", shift left=3, from=1-2, to=1-1]
	\arrow["{s_1}"{description}, shift left=3, from=1-2, to=1-3]
	\arrow["{s_0}"{description}, shift right=3, from=1-2, to=1-3]
	\arrow["{d_2}"', shift right=6, from=1-3, to=1-2]
	\arrow["{d_1}"{description}, from=1-3, to=1-2]
	\arrow["{d_0}", shift left=6, from=1-3, to=1-2]
\end{tikzcd}\]
\begin{proposition}{\cite[Section 4]{GKW}}\label{proposition:GKW}
A unary operadic category is equivalently a pair \hbox{$\bO=(X,\C)$} of a simplicial set $X$ and a category $\C$, such that 
$$ \dec X \cong \N\C.$$
\end{proposition}


\section{The nerve and décalage for 2-categories}\label{section:Nerve_and_dec}
We recall some standard 2-categorical constructions. By a 2-category we always mean a strict 2-category, i.e.~a structure with objects (0-cells), morphisms (1-cells), and morphisms between morphisms (2-cells), with all possible compositions being strictly associative and unital. From now on $\Cat$ denotes the 1-category of 2-categories. It is a monoidal category with the cartesian product. When needed, we regard $(\Cat,\x,\1)$ as an equivalent \textit{strict} monoidal category. In~general, we wish to ignore the coherence isomorphisms, as well as all size issues. 

The 2-categorical nerve was explicitly introduced by J.~Duskin in \cite{Dusk}. It agrees with the classical nerve on 1-categories, regarded as 2-categories with only identity 2-cells, hence we keep the notation $\N$ from the previous section.
\begin{definition}\label{definition:duskin_nerve}
For a 2-category $\C$, its nerve $\N\C$ is a simplicial set defined as $$\N\C_n= \TwoCatNL(\texttt{n},\C),
$$
where $\texttt{n}$ is the ordinal $\texttt{n}=[0 < 1 < \cdots < n]$, viewed as a 2-category, and~$\TwoCatNL$ is the category of 2-categories and lax functors, which strictly preserve identities.
\end{definition}
This gives a fully faithful functor $\N\colon \TwoCat \to \sSet$, whose image consists of is 3-coskeletal simplicial sets. A general $n$-simplex $\sigma$ of $\C$ is a collection of objects $x_0,\ldots, x_n$ of $\C$, maps $f_{ij}\colon x_i \to x_j,$ for $ 0\leq i<j\leq n,$ and lax triangles 
$$ f_{jk}\circ f_{ij}\xRightarrow{\alpha_{ijk}}f_{ik}, \text{ for } 0\leq i<j<k\leq n, $$ such that the following diagram commutes for all $0\leq i<j<k<l\leq n$.
 \[\begin{tikzcd}
	f_{kl}\circ f_{jk}\circ f_{ij} & f_{jl}\circ f_{ij} \\
	f_{kl}\circ f_{ik} & f_{il}
	\arrow["{\alpha_{jkl}\circ f_{ij}}",Rightarrow, from=1-1, to=1-2]
	\arrow["{f_{kl}\circ\alpha_{ijk}}"', Rightarrow,from=1-1, to=2-1]
	\arrow["\alpha_{ijl}", Rightarrow, from=1-2, to=2-2]
	\arrow["{\alpha_{ikl}}", Rightarrow,from=2-1, to=2-2]
\end{tikzcd}\]
Alternatively the nerve can be characterised by strict 2-functors from path 2-categories (cf.~\cite[\href{https://kerodon.net/tag/00B9}{Tag 00B9}]{K}). An $n$-simplex $\sigma$ in $\C$ is thus given by a collection of objects $x_0,\ldots, x_n$ of $\C$, maps $f_S\colon x_i \to x_j$ for every path $S=[i=p_0<\cdots <p_k=j],$ and \hbox{2-cells} $\alpha_{ST}\colon f_S \Rightarrow f_T$ for $T\subseteq S$, an~endpoint preserving (strict) inclusion. It is required that $ f_{S''}\circ f_{S'}= f_{S'\cup S''}$, whenever the paths $S'$ and $S''$ can be concatenated, and that 
\begin{center}
$\alpha_{ST}\circ\alpha_{RS}=\alpha_{TR}$ for $T\subseteq S \subseteq R$.
\end{center}
The first few terms of $\NC$ are the following.
\begin{itemize}
    \item $\NC_0$ is the set of objects of $\C$,
    \item $\NC_1$ is the set of maps $x_0 \xrightarrow{f_{01}}x_1$ of $\C$,
    \item $\NC_2$ is the set of lax triangles
 \[\begin{tikzcd}
	x_0 && x_1\\
	& x_2
	\arrow[""{name=0, anchor=center, inner sep=0},"{f_{01}}", from=1-1, to=1-3]
	\arrow["{f_{02}}"', from=1-1, to=2-2]
	\arrow["{f_{12}}", from=1-3, to=2-2]
	\arrow["\alpha_{012}"{description},draw=none, from=0, to=2-2]
	\end{tikzcd}\]
 in $\C$, where $f_{12} \circ f_{01}\xRightarrow{\alpha_{012}} f_{02}$ is a 2-cell of $\C$,
 \item $\NC_3$ is a set of 3-simplices in $\C$, that is, quadruples $\sigma=(\alpha_{123},\alpha_{023},\alpha_{013},\alpha_{012})$ of \hbox{2-cells} $ f_{jk}\circ f_{ij}\xRightarrow{\alpha_{ijk}}f_{ik},$ such that $\alpha_{013}\circ ({\alpha_{123}\Box f_{01}}) = \alpha_{023}\circ ({f_{23}\Box\alpha_{012}})$ where `$\Box$' is the horizontal composition of $\C$.
\end{itemize} 
We interpret $\sigma$ as a filled 3-simplex obtained by gluing two squares along their boundaries, whose interiors equal:
\[\begin{tikzcd}
	{x_0} & {x_1} \\
	{x_3} & {x_2}
	\arrow["{f_{01}}", from=1-1, to=1-2]
	\arrow[""{name=0, anchor=center, inner sep=0}, "{f_{03}}"', from=1-1, to=2-1]
	\arrow[""{name=1, anchor=center, inner sep=0},  from=1-2, to=2-1]
	\arrow["{f_{12}}", from=1-2, to=2-2]
	\arrow["{f_{23}}", from=2-2, to=2-1]
	\arrow["\alpha_{013}"{description},pos=0.35, draw=none, from=1-1, to=1]
	\arrow["\alpha_{123}"{description},pos=0.4, draw=none, from=2-2, to=1]
\end{tikzcd}
=
\begin{tikzcd}
	{x_0} & {x_1} \\
	{x_3} & {x_2}.
	\arrow["{f_{01}}", from=1-1, to=1-2]
	\arrow[""{name=0, anchor=center, inner sep=0}, "{f_{03}}"', from=1-1, to=2-1]
	\arrow[""{name=1, anchor=center, inner sep=0}, from=1-1, to=2-2]
	\arrow["{f_{12}}", from=1-2, to=2-2]
	\arrow["{f_{23}}", from=2-2, to=2-1]
	\arrow["\alpha_{012}"{description},pos=0.35, draw=none, from=1-2, to=1]
	\arrow["\alpha_{023}"{description},pos=0.4, draw=none, from=2-1, to=1]
\end{tikzcd}
\]
The face maps are clear from the graphical description, $d_i$ deletes the vertex~$x_i$. For a 3-simplex this means $\sigma=(d_0\sigma,d_1\sigma,d_2\sigma,d_3\sigma)$. The degeneracy maps~are:
\begin{itemize}
    \item for an object $x$ in $\C$, $s_0x=\uu_x$,
    \item for a map $f\colon x \to y$ in $\C$, $s_0 f = (f\circ \uu_x \xRightarrow{\uu_f} f)$, and $s_1f = (\uu_y\circ f \xRightarrow{\uu_f} f)$,
    \item for a lax triangle $f_{12}\circ f_{01}\xRightarrow{\alpha_{012}} f_{02},$ of $\C$, 

    $$s_0\alpha=(\alpha,\alpha,\uu,\uu), \, s_1\alpha=(\uu,\alpha,\alpha,\uu), \, s_2\alpha=(\uu,\uu,\alpha,\alpha).$$
\end{itemize}

Next, we study how the upper décalage (Definition~\ref{definition:dec}) interacts with the 2-categorical nerve and we introduce a corresponding construction~`$\D$' for 2-categories.

\begin{definition}
        Let $\C$ be a 2-category and $x$ an object of $\C$. The lax slice $\C/x$ is the defined to be the 2-category with:
    \begin{itemize}
        \item Objects are maps $y\xrightarrow{f}x$ of $\C$ with codomain $x$.
        \item For objects $z\xrightarrow{h} x$ and $y\xrightarrow{g} x$, a map $h\to g$ of $\C/x$ is a pair $(f,\alpha)$, where $ z\xrightarrow{f} y$ is a map of~$\C$ and $g\circ f\xRightarrow{\alpha} h $ is a 2-cell of $\C$.  The pair $(f,\alpha)$ is drawn as:
\[\begin{tikzcd}
	z && y. \\
	& x
	\arrow[""{name=0, anchor=center, inner sep=0},"f", from=1-1, to=1-3]
	\arrow["h"', from=1-1, to=2-2]
	\arrow["g", from=1-3, to=2-2]
	\arrow["\alpha"{description},draw=none, from=0, to=2-2]
\end{tikzcd}\]
        
        \item A 2-cell $\gamma$ of $\C/x$,
        \[\begin{tikzcd}
	h & g
	\arrow[""{name=0, anchor=center, inner sep=0}, "{(f',\alpha')}", curve={height=-12pt}, from=1-1, to=1-2]
	\arrow[""{name=1, anchor=center, inner sep=0}, "{(f'',\alpha'')}"', curve={height=12pt}, from=1-1, to=1-2]
	\arrow["\gamma", shorten <=3pt, shorten >=3pt, Rightarrow, from=0, to=1]
\end{tikzcd}\]
       is a 2-cell $f'\xRightarrow{\gamma} f''$ of $\C$, such that $ \beta\circ (\uu_g \Box \gamma)=\alpha$.
    \end{itemize}
\end{definition}
More details can be found for instance in \cite[Definition~7.1]{JY}, but note the opposite orientation of the triangle interior.  
\begin{definition}
 Let $\C$ be a 2-category and $x$ an object of $\C$. The lax coslice $x/\C$ is the 2-category $$x/\C=(\C^{\op}/x)^{\op},$$ where $(-)^{\op}$ inverts only 1-cells of the 2-category.
\end{definition}

\begin{definition}\label{definition:D}
    For a 2-category $\C$ we write
$$\DC:=\displaystyle\coprod_{x\in \C_0} \C/x.$$
\end{definition}
\begin{example}
Let $\M$ be a strict monoidal category, $\B\M$ its associated 1-object 2-category, and let $U\C$ denote the 1-category obtained from a 2-category $\C$ by forgetting the 2-cells. Recall the category $\M_\triangle$ of Example~\ref{example:O_M}. We have $$\M_\triangle\cong U\D\B\M.$$
\end{example}

\begin{proposition}\label{proposition:dec_and_nerve}
    For a 2-category $\C$,
    $\dec \NC \cong \NDC.$
\end{proposition}
\proof
Since the simplicial set $\NDC$ is 3-coskeletal, it is enough to verify  bijections $$\NDC_k \cong \NC_{k+1}$$ for $k\leq 3$. We comment only on the case $k=3$.

First, a $3$-simplex $\sigma$ of $\DC$ is determined by objects $f_0,\ldots, f_3$ of $\DC$, that is, by maps $x_i\xrightarrow{f_i} y_i$ of $\C$. Further, for every path $S=[i=p_0<\cdots <p_k=j]$, there is a map $f_S=(t_S,\alpha_S)\colon f_i \to f_j$ in $\DC$, that is a lax triangle
 \[\begin{tikzcd}
	x_i && x_j,\\
	& y_i=y_j
	\arrow[""{name=0, anchor=center, inner sep=0},"{t_S}", from=1-1, to=1-3]
	\arrow["{f_i}"', from=1-1, to=2-2]
	\arrow["{f_j}", from=1-3, to=2-2]
	\arrow["\alpha_S"{description},draw=none, from=0, to=2-2]
	\end{tikzcd}\]
 and for every subsequence $T\subseteq S$ there is a 2-cell $\alpha_{ST}\colon f_S\Rightarrow f_T$ in $\DC$, that is, a 2-cell $\alpha_{ST}\colon t_S\Rightarrow t_T$ in $\C$, such that $\alpha_T\circ (f_j\square \alpha_{ST})=\alpha_S$. 
On the other hand, a 4-simplex $\delta$ of $\C$ is given by objects $z_0,\ldots,z_4$, maps $g_S\colon z_i \to z_j$, and 2-cells $\beta_{ST}\colon g_S \to g_T $ of $\C$.

Given $\delta$ we construct a 3-simplex $\sigma$ of $\DC$ as follows. Define
$ f_i := g_{[i<4]}\colon z_i \to z_4$. Let~$S$ be the path $S=[i=p_0<\cdots <p_k=j]$ and denote $\Bar{S}:=S\cup [j<4]$. Since $g_{[j<4]}\circ g_S = g_{\Bar{S}}$, we define $f_S:=\beta_{[i<4]\Bar{S}}\colon g_{[j<4]}\circ g_S \to g_{[i<4]} $ and $\alpha_{ST}:=\beta_{ST}$.

Conversely, given a 3-simplex $\sigma$ of $\DC$, we construct $\delta$ as follows. For $i\leq 3$ define $z_i:=x_i$ and $z_4:=y_0=\cdots =y_4$. Further, $g_S:=t_S$ and $\beta_{ST}:=\alpha_{ST}$. These two constructions are inverse and the bijection preserves the simplicial structure.
\endproof

\section{Unary operadic 2-categories}\label{section:UnOpBic}
The following definition was motivated by the characterization of unary operadic categories of \cite{GKW}, which we recalled in Proposition~\ref{proposition:GKW}. We adapt the characterization to the 2-categorical setting and we use it as a starting point for our theory.

 \begin{definition}\label{definition:unary_operadic_2-category_A}
     A unary operadic 2-category is a pair $\bO=(X,\C)$ of a 2-category $\C$ and a simplicial set $X$, such that the upper décalage of $X$ is the nerve of $\C$:
     $$\dec X\cong \NC.$$
    An operadic functor between two unary operadic 2-categories is a morphism of their underlying simplicial sets. 
We~denote the category of unary operadic 2-categories by $\TwoUnOpCat$.
\end{definition}
We unpack the definition in terms of the underlying 2-category $\C$, the set $X_0$, and additional top face and degeneracy maps on $\NC$, which we denote by $\phi$ and $u$, respectively.

\begin{proposition}\label{proposition:unary_operadic_2-category_B}
A unary operadic 2-category is a 2-category $\C$ together with the following~data:
\begin{enumerate}
        \item \label{item:1}for every object $a$ of $\C$, there is $\phi(a)\in \pi_0(\C)$,
        \item \label{item:2}for every arrow $x\xrightarrow{q}y$ in $\C$, $\phi(q)$ is an object of $\C$,
        \item \label{item:3}for every lax triangle $g\circ f \xRightarrow{\alpha} h$ in $\C$, $\phi(h)\xrightarrow{\phi(\alpha)}\phi(g)$ is a map in $\C$,
        \item \label{item:4}for a 3-simplex $\sigma_{0123}$ of $\C$,  $\phi(\alpha_{123})\circ \phi(\alpha_{013}) \xRightarrow{\phi(\sigma_{0123})}\phi(\alpha_{023})$ is a lax triangle in $\C$, 
        \item \label{item:5}for every $c \in \pi_0(\C)$, $u_c$ is an object of $\C$,
    \item \label{item:6}for every object $a\in \C$, $a \xrightarrow{u_a} u_{\pi(a)}$ is a map in $\C$,
    \item \label{item:7}for every map $x\xrightarrow{q} y$ of $\C$, there is a 2-cell $u_y\circ q \xRightarrow{u_q} u_x$, and
\item \label{item:8}for every lax triangle $ g\circ f\xRightarrow{\alpha}h$, there is a 3-simplex $u_\alpha=(u_g,u_h,u_f,\alpha)$ in $\C$. 
\end{enumerate}
This data are subject to the following equations. For every $c, a, q, \alpha$, and $\sigma_{0123}$ as above,  
\begin{enumerate} 
\setcounter{enumi}{8}
\item \label{item:9} $ \pi(\phi (q))=\phi(x)$,
    \item\label{item:10} $\pi(u_c)=c$,
    \item\label{item:11}
    $\phi(\uu_a)=u_{\phi(a)}$,
    \item\label{item:12}  $\phi\big( \uu_y\circ q \xRightarrow{\uu_q}q\big)=u_{\phi(q)}$ and $\phi\big(q\circ \uu_x \xRightarrow{\uu_q}q\big)=\uu_{\phi(q)}$,
    
    \item\label{item:13}
   $\phi(\uu,\uu,\alpha,\alpha)=u_{\phi(\alpha)}$,
   $\phi(\uu,\alpha,\alpha,\uu)=\uu_{\phi(g)}\circ\phi(\alpha)\xRightarrow{\uu_{\phi(\alpha)}} \phi(\alpha),$ and
   $\phi(\alpha,\alpha,\uu,\uu)=\phi(\alpha)\circ \uu_{\phi(h)} \xRightarrow{\uu_{\phi(\alpha)}}\phi(\alpha),$
   
    \item\label{item:14} $u_{u_c}=\uu_{u_c}$, $u_{u_a}=\uu_{u_a}$, and $u_{u_q}=\uu_{u_q}$, 
    \item\label{item:15} $\phi\circ u = \uu$, more precisely
    $\phi(u_c)=c$,
    $\phi(u_a)=a$,
     $\phi(u_q)=q$, and
    $\phi(u_\alpha)=\alpha$,
    \item\label{item:16} $\phi(\phi(\sigma_{0123}))=\phi(\alpha_{012})$, $\phi(\phi(\alpha))=\phi(f)$, and $\phi(\phi(q))=\phi(x)$,
    \item\label{item:17} $u_{\uu_x}=\uu_{u_x}$, $u_{s_0 q}=s_0 u_q$, and $u_{s_1 q}=s_1 u_q$. 
    
\end{enumerate}
\end{proposition}
\begin{remark}
These conditions are not minimal. For example the condition $$\phi(\phi(q\colon x\to y))=\phi(x)$$ can be derived from $\phi(\phi(\alpha))=\phi(f)$ with $\alpha = s_0 q$ as follows:
$$\phi\phi q = \pi u\phi\phi q=\pi\phi s_0\phi q=\pi\phi\phi s_0 q =\pi\phi d_2 s_0 q = \pi \phi s_0 d_1 q = \pi u \phi d_1 q = \phi d_1 q = \phi x.
$$
\end{remark}
Before we prove the proposition, we make two observations. Let $\bO=(\C,X)$ be a unary operadic 2-category.
\begin{lemma} There is a bijection 
    $X_0\cong \pi_0(\C)$, such that the following diagram commutes.
   \[\begin{tikzcd}
	X_0 & \pi_0(\C) \\
	X_1 &  \C_0
	\arrow["\cong", from=1-1, to=1-2]
	\arrow["\cong", from=2-1, to=2-2]
	\arrow["d_0", from=2-1, to=1-1]
	\arrow["\pi", from=2-2, to=1-2]
 \end{tikzcd}\] 
\end{lemma}
\proof
This follows from the fact that for any simplicial set X, the set $X_0$ is a split coequalizer of $d_0,d_1\colon X_2 \rightrightarrows X_1$ and $\pi_0(\C)$ is a coequalizer of $X_2=\C_1 \rightrightarrows \C_0=X_1$.
\endproof
\begin{lemma}
The simplicial set $X$ is 4-coskeletal.
\end{lemma}
\proof
This follows from the defining isomorphism $\dec X\cong \NC$ and the fact that the nerve~$\NC$ is 3-coskeletal. Let the isomorphism consist of bijections
$X_{p+1}\xrightarrow{\iota_p} \NC_p$, for $p\geq 0,$
with $$\iota_p d_i=d_i \iota_{p+1},$$ for $0\leq i \leq p+1$. 
Any boundary $\delta\Delta [k]\xrightarrow{\sigma}X$ of a $k$-simplex in~$X$, for $k>4$, restricts to a horn $\sigma_{\hat{k}}$ by forgetting the face $d_k\sigma$. That is, $\sigma_{\hat{k}}$ is a $k$-tuple of compatible $(k-1)$-cells  
$$\sigma_{\hat{k}}=(d_0\sigma,\ldots,d_{k-1}\sigma).$$
Since $d_i\sigma\in X_{k-1}=(\dec X)_{k-2}$, the horn $\sigma_{\hat{k}}$ determines a boundary of a $(k-1)$-simplex in $\dec X$, which has a unique filler $\bar{\sigma} \in (\dec X)_{k-1}=X_k$, with $d_i\sigma=d_i\bar{\sigma}$ for $0\leq i\leq k-1$, which means $d_i\sigma=d_i\bar{\sigma}$ in $X_{k-1}$. We need to prove that also $d_k\sigma=d_k\bar{\sigma}$. Since $\iota$ is a bijection, it~is enough to show $\iota_{k-2}d_k\sigma=\iota_{k-2}d_k\bar{\sigma}\in \NC_{k-2}$. Any element $\alpha\in\NC_{k-2}$ is a $(k-1)$-tuple $\alpha=(d_0\alpha,\ldots,d_{k-2}\alpha)$ of compatible 2-cells, and two elements $\alpha,\beta \in \NC_{k-2}$ equal if, and only if $d_i\alpha=d_i\beta$ for $0 \leq i \leq k-2$. We compute  
$$
d_i\iota_{k-2} d_k\sigma = \iota_{k-3}  d_i d_k \sigma = \iota_{k-3} d_{k-1} d_i \sigma=\iota_{k-3} d_{k-1} d_i \bar{\sigma}=\iota_{k-3} d_i d_{k} \bar{\sigma}=d_i \iota_{k-2} d_{k} \bar{\sigma}, 
$$
which shows that $\bar{\sigma}$ is a unique filler of $\sigma$ in $X_k$.
\endproof 
\begin{proof}[Proof of Proposition~\ref{proposition:unary_operadic_2-category_B}]
 It follows from the two lemmas above that a unary operadic 2-category is determined by a diagram 
\[\begin{tikzcd}
	{\pi_0\C} & {\C_0} & {\C_1} & {\NC_2} & {\NC_3}
	\arrow["u"{description}, from=1-1, to=1-2]
	\arrow["\phi"{description}, shift right=3, from=1-2, to=1-1]
	\arrow["\pi"{description}, shift left=3, from=1-2, to=1-1]
	\arrow["u"{description}, shift left=2, from=1-2, to=1-3]
	\arrow[shift right=2, from=1-2, to=1-3]
	\arrow["\phi"{description}, shift right=5, from=1-3, to=1-2]
	\arrow["{d_0}", shift left=4, from=1-3, to=1-2]
	\arrow[from=1-3, to=1-2]
	\arrow["u"{description}, shift left=2, from=1-3, to=1-4]
	\arrow["\phi"{description}, shift right=5, from=1-4, to=1-3]
	\arrow[""{name=0, anchor=center, inner sep=0}, from=1-4, to=1-3]
	\arrow[""{name=1, anchor=center, inner sep=0}, "{d_0}", shift left=4, from=1-4, to=1-3]
	\arrow["u"{description}, shift left=2, from=1-4, to=1-5]
	\arrow["\phi"{description}, shift right=5, from=1-5, to=1-4]
	\arrow[""{name=2, anchor=center, inner sep=0}, from=1-5, to=1-4]
	\arrow[""{name=3, anchor=center, inner sep=0}, "{d_0}", shift left=4, from=1-5, to=1-4]
	\arrow["\cdots"{marking, allow upside down}, draw=none, from=0, to=1]
	\arrow["\cdots"{marking, allow upside down}, draw=none, from=2, to=3]
\end{tikzcd}\]
 and simplicial identities up to $X_4=(\NC)_3$.
Most of the identities are satisfied since $\C$ is a 2-category. We translate those which are not automatic:\\
\begin{table}[H]
\resizebox{\textwidth}{!}{%
\begin{tabular}{ccccc}
\begin{tabular}{|c|} 
 \hline
 $\C_1\to \pi_0\C$\\
  \hline
 $\pi\phi=\phi d_0$ \eqref{item:9}\\
 $\phi\phi=\phi d_1$ \eqref{item:16}\\
 \hline
\end{tabular}
&
\begin{tabular}{|c|} 
 \hline
 $\NC_2\to \C_0$\\
  \hline
 $d_0\phi=\phi d_0$ \eqref{item:3}\\
 $d_1\phi=\phi d_1$ \eqref{item:3}\\
 $\phi\phi=\phi d_2$ \eqref{item:16}\\
 \hline
\end{tabular}
&
\begin{tabular}{|c|} 
 \hline
 $\NC_3\to \C_1$\\
 \hline
 $d_0\phi=\phi d_0$ \eqref{item:4}\\
 $d_1\phi=\phi d_1$ \eqref{item:4}\\
 $d_2\phi=\phi d_2$ \eqref{item:4}\\
 $\phi\phi=\phi d_3$ \eqref{item:16}\\
 \hline
\end{tabular}
&
\begin{tabular}{|c|} 
 \hline
 $\C_1\to \NC_3$\\
 \hline
 $u s_0=s_0u$ \eqref{item:17}\\
 $u s_1=s_1u$ \eqref{item:17}\\
 $u u=s_2u$ \eqref{item:14}\\
 \hline
\end{tabular}
&
\begin{tabular}{|c|} 
 \hline
 $\C_0\to \NC_2$\\
 \hline
 $u s_0=s_0u$ \eqref{item:17}\\
 $u u=s_1u$ \eqref{item:14}\\
 \hline
\end{tabular}
\\
\begin{tabular}{|c|} 
 \hline
 $\pi_0\C\to \C_1$\\
 \hline
 $uu=s_0u$ \eqref{item:14}\\
 \hline
\end{tabular}
&
\begin{tabular}{|c|} 
 \hline
 $\pi_0\C\to \pi_0\C$\\
 \hline
 $\pi u=1$ \eqref{item:10}\\
 $\phi u=1$ \eqref{item:15}\\
 \hline
\end{tabular}
&
\begin{tabular}{|c|} 
 \hline
 $\C_0\to \C_0$\\
 \hline
 $\phi s_0=u\phi$ \eqref{item:11}\\
 $d_0u=u \pi$ \eqref{item:6}\\
 $d_1u=1$ \eqref{item:6}\\
 $\phi u=1$ \eqref{item:15}\\
 \hline
\end{tabular}
&
\begin{tabular}{|c|} 
 \hline
 $\C_1\to \C_1$\\
 \hline
 $\phi s_0=s_0\phi$ \eqref{item:12}\\
 $\phi s_1=u\phi$ \eqref{item:12}\\
 $d_0u=u d_0$ \eqref{item:7}\\
 $d_1u=u d_1$ \eqref{item:7}\\
 $d_2u=1$ \eqref{item:7}\\
 $\phi u=1$ \eqref{item:15}\\
 \hline
\end{tabular}
&
\begin{tabular}{|c|} 
 \hline
$\NC_2\to \NC_2$\\
 \hline
 $\phi s_0=s_0\phi$  \eqref{item:13} \\
  $\phi s_1=s_1\phi$  \eqref{item:13}  \\
    $\phi s_2=u \phi$ \eqref{item:13}\\
    $d_0u=u d_0$ \eqref{item:8}\\
 $d_1u=u d_1$ \eqref{item:8}\\
 $d_2u=u d_2$ \eqref{item:8}\\
 $d_3u=1$ \eqref{item:8}\\
 $\phi u=1$ \eqref{item:15}\\
 \hline
\end{tabular}
\end{tabular}%
}
\end{table}

 These equations compare directly (cf.~numbers in brackets) to the data and conditions of Proposition~\ref{proposition:unary_operadic_2-category_B}, which finishes the proof.
\end{proof}
Our next goal is to reformulate Definition~\ref{definition:unary_operadic_2-category_A}, so it resembles the original definition of unary operadic categories of \cite{Batanin_Markl:blob}, see Definition~\ref{definition:unary_operadic_category}. We find analogues to fiber functors and local terminal objects.
\begin{definition}
    An object $v$ in a 2-category $\C$ is lali-terminal, 
    if for every object~$x$ of $\C$, the category $\C(x,v)$ has a terminal object. We require that the terminal object of $\C(v,v)$ is the identity on $v$.
    An object $v$ is local lali-terminal if it is lali-terminal in its component. 
\end{definition}
The prefix \textit{lali} stands for `left adjoint left inverse'. This terminology appears in \cite{Step}, however the condition on the terminal object of $\C(v,v)$ being the identity is not required there.
\begin{lemma}\label{lemma:lali}
    An object $v$ in a 2-category $\C$ is lali-terminal if and only if for every object $x\in \C$, there is a chosen map $x \xrightarrow{v_x} v$, for every map $x\xrightarrow{f} y$ of $\C$, there is a chosen 2-cell $v_y\circ f \xRightarrow{v_f} v_x$, and for every lax triangle $ g\circ f\xRightarrow{\alpha}h$, there is a 3-simplex $v_\alpha=(v_g,v_h,v_f,\alpha)$ in $\C$, such that $v_v=\uu_v$ and $v_{v_x}=\uu_{v_x}$. 
\end{lemma}
\proof
Let $v$ be lali-terminal in $\C$. For an object $x$, the map $v_x$ is the terminal object of $\C(x,v)$ and for $f\colon x \to y$, the 2-cell $v_f$ is the terminal map $v_y\circ f\Rightarrow v_x$ of $\C(x,v)$. For a 2-cell $g\circ f \xRightarrow{\alpha} h$, there are two 2-cells $v_z\circ g\circ f \Rightarrow v_x$, namely $v_h \circ ( v_x\cdot\alpha )$ and $v_f \circ (v_g\cdot f )$. Since $v_x$ is terminal, the two composites equal, hence they form a 3-simplex $v_\alpha=(v_g,v_h,v_f,\alpha)$. By definition,  $v_v=\uu_v$ and it follows that $v_{v_x}=\uu_{v_x}$.

Conversely, we will show that the choice of $v_x,v_f$, and $v_\alpha$, with $v_v=\uu_v$ and $v_{v_x}=\uu_{v_x}$, makes $v$ lali-terminal. Let there be two 2-cells $\alpha, \beta\colon f\Rightarrow v_x$ in $\C(x,v)$. We interpret them as lax triangles $f\circ \uu_x \Rightarrow v_x$. We have their 3-simplices $v_\alpha$ and $v_\beta$, which give $$\alpha=v_{v_x} \circ (v_v\cdot \alpha)= v_{\uu_x} \circ (v_f\cdot \uu_x )=v_{v_x} \circ (v_{v}\cdot \beta)=\beta.$$
\endproof

\begin{proposition}\label{proposition:unary_operadic_2-category_C}
    A unary operadic 2-category is a 2-category $\C$ together with a normalised lax functor $\D\C\xrightarrow{\phi} \C$ and a set $\mathscr{U}$ of chosen local lali-terminal objects $u_c$ in each connected component $c$ of $\C$, such that the following conditions hold.
 \begin{itemize}
    \item For any object $x$ of $\C$, $\phi(\uu_x) \in \mathscr{U},$
    \item for a map $f\colon x\to y$ of $\C$,
    $\phi(f,\uu_f)=u_{\phi(f)},$
    where  $(f,\uu_f)$ is the 2-cell 
    $ \uu_y\circ f \xRightarrow{\uu_f}f$,
    \item for $u_c \in \mathscr{U}$, $\phi_{u_c}\colon \C/u_c\to \C$ is the domain functor,
    \item for a 3-simplex $\sigma_{0123}$ in $\C$, $\phi(\phi(\sigma_{0123}))=\phi(\alpha_{012})$.
     
 \end{itemize}   
\end{proposition}
\proof
Elementary, using the  characterisation of Proposition~\ref{proposition:unary_operadic_2-category_B} and Lemma~\ref{lemma:lali}.
\endproof
\begin{example}
If we consider 1-categories as 2-categories with only identity 2-cells, then    every unital unary operadic category of \cite{Batanin_Markl:blob} 
    is a unary operadic 2-category. We shall call them unary operadic 1-categories.
\end{example}
\begin{example}\label{example:monoidal_categories}
By Proposition~\ref{proposition:dec_and_nerve}, for any 2-category $\C$, the pair $\bC=(\NC,\DC)$ is a unary operadic 2-category. In particular, any strict monoidal category $\V$ can be viewed as a one-object 2-category $\B\V$. Hence, $(\N\B\V,\D\B\V)$ is a unary operadic 2-category which we denote by~$\V_\triangle$. This is in a slight conflict with the notation of Example~\ref{example:O_M}, however we stick to this new meaning. Let us make the structure explicit. The 2-category~$\D\B\V$ consists of the following data.
     \begin{itemize}
         \item Objects are the objects of $\V$.
         \item Morphisms $x\to y$ are pairs $(a,f)$ such that $y \ox a\xrightarrow{f}  x$ is a map in $\V$.
         \item 2-cells
         $(a,f)\xRightarrow{\delta}(b,g)$ are maps $\delta\colon a \to b $ of $\V$, such that the triangle
\[\begin{tikzcd}
	y\ox a && y\ox b\\
	& x
	\arrow["{y\ox \delta}", from=1-1, to=1-3]
	\arrow["{f}"', from=1-1, to=2-2]
	\arrow["{g}", from=1-3, to=2-2]
	\end{tikzcd}\]
 commutes.
          \end{itemize}
          The additional structure maps $u$ and $\phi$ on $\N\D\B\V$ are obtained from the simplicial set $\N\B\V$ and the isomorphism $\dec \N \B\V \cong \N\D\B\V$. Namely, $\N\B\V$ consists of:
\begin{itemize}
    \item $(\N\B\V)_0=\{\cdot\}$, 
    \item $(\N\B\V)_1$ is the set of segments
    $$\cdot\xrightarrow{x} \cdot,\, x \in \V,$$
    \item $(\N\B\V)_2$ is the set of triangles
\[\begin{tikzcd}
	\cdot && \cdot\\
	& \cdot
	\arrow[""{name=0, anchor=center, inner sep=0}, "{a_{01}}", from=1-1, to=1-3]
	\arrow["{x_0}"', from=1-1, to=2-2]
	\arrow["{x_1}", from=1-3, to=2-2]
	\arrow["f_{01}"{description}, draw=none, from=2-2, to=0]
\end{tikzcd}\]
    with a map $ x_1 \ox a_{01}\xrightarrow{f_{01}} x_0 $ in $\V$,
     \item $(\N\B\V)_3$ is the set of 3-simplices 

     \[\begin{tikzcd}
	\cdot & \cdot \\
	\cdot & \cdot
	\arrow["{a_{01}}", from=1-1, to=1-2]
	\arrow["{x_0}"', from=1-1, to=2-1]
	\arrow[""{name=1, anchor=center, inner sep=0},"{x_1}"{description},  from=1-2, to=2-1]
	\arrow["{a_{12}}", from=1-2, to=2-2]
	\arrow["{x_{2}}", from=2-2, to=2-1]
	\arrow["f_{01}"{description},pos=0.25, draw=none, from=1-1, to=1]
	\arrow["f_{12}"{description},pos=0.25, draw=none, from=2-2, to=1]
\end{tikzcd}
=
\begin{tikzcd}
\cdot & \cdot \\
	\cdot & \cdot,
	\arrow["{a_{01}}", from=1-1, to=1-2]
	\arrow["{x_0}"', from=1-1, to=2-1]
	\arrow[""{name=1, anchor=center, inner sep=0},"{a_{02}}"{description},  from=1-1, to=2-2]
	\arrow["{a_{12}}", from=1-2, to=2-2]
	\arrow["{x_{2}}", from=2-2, to=2-1]
	\arrow["\alpha_{012}"{description},pos=0.25, draw=none, from=1-2, to=1]
	\arrow["f_{02}"{description},pos=0.25, draw=none, from=2-1, to=1]
\end{tikzcd}
\]
\item (we do not specify $(\N\B\V)_4$, instead we describe $(\N\D\B\V)_3$ directly below.)  
\end{itemize}
The simplicial set $\N\D\B\V$ consists of:
\begin{itemize}
    \item $(\N\D\B\V)_0=\V_0$, 
    \item $(\N\D\B\V)_1$ is the set of segments
    $$x_0\xrightarrow{(a_{01},f_{01})} x_{1},$$
    where $x_1 \ox a_{01}\xrightarrow{f_{01}}  x_0$ is a map in $\V$,
    \item $(\N\D\B\V)_2$ is the set of triangles
\[\begin{tikzcd}
	x_0 && x_1\\
	& x_2
	\arrow["{(a_{01},f_{01})}", from=1-1, to=1-3]
	\arrow[""{name=0, anchor=center, inner sep=0}, "{(a_{02},f_{02})}"', from=1-1, to=2-2]
	\arrow["{(a_{12},f_{12})}", from=1-3, to=2-2]
	\arrow["\alpha_{012}"{description},pos=0.4, draw=none, from=0, to=1-3]
\end{tikzcd}\]
    with a map $ a_{12}\ox a_{01} \xrightarrow{\alpha_{012}} a_{02} $ in $\V$, such that the diagram
    \[\begin{tikzcd}
	x_{2}\ox a_{12}\ox a_{01} & x_2\ox a_{02} \\
	x_1\ox a_{01} & x_0
	\arrow["{x_{2}\ox \alpha_{012}}", from=1-1, to=1-2]
	\arrow["{f_{12}\ox a_{01}}"', from=1-1, to=2-1]
	\arrow["f_{02}", from=1-2, to=2-2]
	\arrow["f_{01}", from=2-1, to=2-2]
\end{tikzcd}\]
commutes. We denote the triangle by $\Delta_{012}=[f_{12},f_{02},f_{01},\alpha_{012}]$, the square brackets indicate the commuting square.
    \item $(\N\D\B\V)_3$ is the set of quadruples $\sigma_{0123}=(\Delta_{123},\Delta_{023},\Delta_{013},\Delta_{012})$ of triangles, that is, a collection of maps $$x_q\ox a_{pq}\xrightarrow{f_{pq}} x_p\text{ and } a_{jk}\ox a_{ij}\xrightarrow{\alpha_{ijk}}a_{ik}$$
in $\V$, such that the following cube commutes. The labeling of each edge is determined by its source and target.
\[\begin{tikzcd}
	{x_3\ox a_{23}\ox a_{12}\ox a_{01}} && {x_2\ox a_{12}\ox a_{01}} \\
  & {x_3\ox a_{13}\ox a_{01}} && {x_1\ox a_{01}} \\
	{x_3\ox a_{23}\ox a_{02}} && {x_2\ox a_{02}} \\
	& {x_3\ox a_{03}} && {x_0}
	\arrow[from=1-1, to=1-3]
	\arrow[from=1-1, to=2-2]
	\arrow[from=1-1, to=3-1]
	\arrow[from=1-3, to=2-4]
	\arrow[dashed, from=1-3, to=3-3]
	\arrow[from=2-2, to=2-4]
	\arrow[from=2-2, to=4-2]
	\arrow[from=2-4, to=4-4]
	\arrow[dashed, from=3-1, to=3-3]
	\arrow[from=3-1, to=4-2]
	\arrow[dashed, from=3-3, to=4-4]
	\arrow[from=4-2, to=4-4]
\end{tikzcd}\]
\end{itemize}

Note that $\D\B\V$ has only one connected component $c_\1$, since for any $x \in \V$ there is a map $x\xrightarrow{(x,\uu_x)}\1$ into the monoidal unit $\1$ of $\V$.
Hence, the nerve $\N\D\B\V$ is indeed extended by the set $\pi_0(\N\D\B\V)= \{c_\1\}\cong (\N\B\V)_0$. The maps $u$ and $\phi$ on $\N\D\B\V$ are:
\begin{itemize}
    \item on connected components: $u_{c_\1}=\1$,
    \item on 0-simplices: $u_x= (x,\uu_x)$ and $\phi(x)=c_\1$,
    \item on 1-simplices: $u_{(a_{01},f_{01})} = [f_{01},f_{01},\uu_{x_0},\uu_{x_1}]$ and $\phi(a_{01},f_{01}) = a_{01}$,
    \item on 2-simplices: $u_{\Delta_{012}}=(\Delta_{012}, u_{(a_{01},f_{01})},u_{(a_{02},f_{02})},u_{(a_{12},f_{12})})$ and $\phi(\Delta_{012}) = (a_{01},\alpha_{012})$,
    \item on 3-simplices: $\phi(\sigma_{0123})=\Delta_{012}$.
\end{itemize}
The 2-category $\D\B\V$ is a 2-categorical version of the para-construction, mentioned in Example~\ref{example:O_M}, which might be of interest.
For an ordinary monoid $M$ in $\Set$, viewed as a discrete monoidal category, $\N\D\B\M$ recovers the classical bar complex \hbox{$(\N\D\B M)_n=M^{\x n}$.} 
\end{example}
\section{Operads of unary operadic 2-categories}
Before we introduce $\bO$-operads of an operadic 2-category $\bO$, we define their 1-connected variant. An operad $\tP$ is 1-connected if the part $\tP_{u_c}$ over chosen local (lali-)terminal object~$u_c$ contains only the identity operation $e_c$. The reason is, that 1-connected \hbox{$\bO$-operads} arise from the simplicial viewpoint naturally. We then distill an explicit definition of an operad and weaken the \hbox{1-connected} condition.

\begin{definition}\label{definition:operad}
Let $\V$ be a strict monoidal category and $\V_\triangle$ the unary operadic 2-category of~Example~\ref{example:monoidal_categories}.
A 1-connected $\bO$-operad $\tP$ with values in $\V$ is an operadic functor \hbox{$\tP\colon\bO\to\V_\triangle$}.
\end{definition}
For simplicity let $\V=\Cat$, the category of small categories (considered as an equivalent strict monoidal category) and let $\bO=(X,\C)$. A categorical $\bO$-operad $\tP$ is thus a collection of categories~$\tP_x$, indexed by objects $x \in \C$, such that $\tP_{u_{c}}= \1$, the terminal category, for every $c\in \pi_0(\C)$. We denote the unique object of $\tP_{u_c}$ by $e_c$. The collection is further equipped with functors $$\tP_y \x \tP_{\phi(f)} \xrightarrow{\tP_f} \tP_x,$$
for every map $f\colon x \to y$ in $C$
 which satisfy the following associativity and unit laws. We~use the notation $a\cdot_f b$ for $\tP_f(a,b)$.
\begin{itemize}
    \item For any lax triangle $\alpha$ in $C$,
    \[\begin{tikzcd}
	x && y\\
	& z
	\arrow[""{name=0, anchor=center, inner sep=0},"{f}", from=1-1, to=1-3]
	\arrow["{h}"', from=1-1, to=2-2]
	\arrow["{g}", from=1-3, to=2-2]
	\arrow["\alpha"{description},draw=none, from=2-2, to=0]
	\end{tikzcd}\]
    and objects $a\in \tP_z, b\in \tP_{\phi(g)}, c \in \tP_{\phi(f)},$
    $$ (a\cdot_{g} b) \cdot_{f}c = a \cdot_{h} (b \cdot_{\phi(\alpha)} c).$$
    It follows from the axioms of an operadic 2-category that $\phi(f)=\phi(\phi(\alpha))$, hence the expression above is well defined.
    \item For any object $x\in C$ and $a\in \tP_x$,
    $$ a\cdot_{\uu_x}e_{\phi(x)}=a ,$$
    $$ e_{\pi(x)}\cdot_{u_x}a=a .$$
\end{itemize}
From the simplicial viewpoint, the operad multiplication has a more pleasing form:
$$
\tP_{d_0f}\x \tP_{d_2f}\xrightarrow{\tP_f} \tP_{d_1f},
$$
and the laws can be written as
$$
\tP_{d_1\alpha}(\tP_{d_3\alpha} \x \uu)= \tP_{d_2\alpha}(\uu \x \tP_{d_0\alpha})
$$
$$
 \tP_{s_0x}(-,e_{d_1x})=\uu_{\tP_x}= \tP_{s_1x}(e_{d_0x},-).
$$
The higher information of Definition~\ref{definition:operad} is redundant. This motivates the following definition where we weaken the 1-connected condition.
\begin{definition}\label{definition:operad2}
A categorical $\bO$-operad $\tP$ is a structure described under Definition~\ref{definition:operad} with the condition $\tP_{u_c}=\{e_c\}$ replaced by a choice of maps $1\xrightarrow{e_c} \tP_{u_c}$ for every $c\in\pi_0(\C)$, which pick the units.
\end{definition}

\begin{example}
    For unary operadic 1-categories we recover the original definition of a unital $\bO$-operad of \cite{Batanin_Markl:blob}.
\end{example}
\begin{example}
    For a 2-category $\C$ denote $\bC=(\NC,\D\C)$. The $\bC$-operads in $\V$ are lax 2-functors from $\C$ to $\B\V$.
    In particular, for a strict monoidal category $\M$, the \hbox{$\M_\triangle$-operads} in $\V$ are lax monoidal functors $\M\to\V$.
\end{example}

\section{Operadic Grothendieck construction}
For a unary operadic 2-category $\bO$ and a 1-connected categorical $\bO$-operad $\tP$, we introduce its operadic Grothendieck construction, i.e.~a unary operadic 2-category $\int_\bO \tP$ with a projection $\pi_\tP\colon \int_\bO \tP\to \bO$. We make the construction explicit. For non-1-connected operads the construction is obtained simply by weakening the 1-connected condition.

Observe that for any $\C\in \Cat$, the lax coslice category $\C/\Cat$ is again monoidal. We~keep the convention that lax means `from many to one'. The morphisms $F\to G$ of $\C/\Cat$ are pairs $(H,\alpha)$: 
\[\begin{tikzcd}
	 &\C& \\
\mathscr{D} && \mathscr{E},
	\arrow[""{name=0, anchor=center, inner sep=0},"{H}"', from=2-1, to=2-3]
	\arrow["{F}"', from=1-2, to=2-1]
	\arrow["{G}", from=1-2, to=2-3]
	\arrow["\alpha"{description},draw=none, from=0, to=1-2]
	\end{tikzcd}\]
 with $H \circ F\xRightarrow{\alpha} G$.
The codomain projection $\pi\colon \C/\Cat \to \Cat$ is a monoidal functor which induces an operadic functor $\pi_\triangle\colon (\C/\Cat)_\triangle \to \Cat_\triangle$ (for the notation cf.~Example~\ref{example:monoidal_categories}).
\begin{proposition}
The category $\TwoUnOpCat$ has pullbacks.
\end{proposition}
\pf Let $p\colon \bO \to \bP$ and $q\colon \bQ \to \bP$ be two operadic functors, with 
\begin{center}
$\bO=(X_\bO,\C_\bO)$, $\bP=(X_\bP,\C_\bP)$, and $\bQ=(X_\bQ,\C_\bQ)$.
\end{center}
It follows from the properties of the nerve and décalage that the pullback of $p$ and $q$ is $$\bO\x_{\bP} \bQ = (X_\bO\x_{X_\bP} X_\bQ,\C_\bO\x_{\C_\bP} \C_\bQ).\qedhere$$
\begin{definition}\label{definition:G-construction} Let $\bO=(X,\C)$ be a unary operadic 2-category and $\tP$ a categorical \hbox{1-connected} $\bO$-operad. The operadic Grothendieck construction 
      $\int_\bO \tP=(X_\tP,\C_\tP)$ is a pullback of $\tP$ along $\pi_\triangle \colon (\1/\Cat)_\triangle \to \Cat_\triangle$, where $\1$ denotes the terminal category:
\[\begin{tikzcd}
	\int_\bO \tP & (\1/\Cat)_\triangle \\
	\bO & \Cat_\triangle.
	\arrow[from=1-1, to=1-2]
	\arrow["\pi_\tP"',from=1-1, to=2-1]
	\arrow["\lrcorner"{anchor=center, pos=0.125}, draw=none, from=1-1, to=2-2]
	\arrow["\pi_\triangle",from=1-2, to=2-2]
	\arrow["\tP",from=2-1, to=2-2]
\end{tikzcd}\]

\end{definition}

The assignment $\int_\bO$ induces a fully faithful functor $\oper{\bO}{\Cat}\to \TwoUnOpCat/\bO$. Its codomain is a category of operadic functors into $\bO$ and morphisms are commutative triangles of operadic functors.
We make the definition explicit in terms of the 2-category~$\C_\tP$ and additional maps $\phi$ and $u$ on $\N\C_\tP$. Denote by $\psi$ the fiber assignment of $\bO$, by $v$ the unit maps of $\bO$, and by $e_c \in \tP_{v_c}$ the units of the operad $\tP$.
\begin{itemize}
    \item The objects of $\C_\tP$ are pairs $(A\in \C_0,a\in \tP_A)$,
    \item 1-cells $(A,a)\to (B,b)$ are triples $(f, p,\alpha)$, where $f\colon A\to B$ is a map in $\C$, $p\in \tP_{\psi(f)}$, and $b\cdot_f p \xrightarrow{\alpha} a $ is a map in $\tP_A$. The identity 1-cell of an objects $(A,a)$ is the triple $(\uu_A,e_{\phi(A)},\uu_a)$. The composite of 1-cells 
    $$
    (A,a)\xrightarrow{(f', p',\alpha')} (B,b) \xrightarrow{(f'', p'',\alpha'')} (C,c)
    $$
    is defined by
    $$
    (f'', p'',\alpha'')\circ 
    (f', p',\alpha') 
    = \big(f''\circ f', p''\cdot_{\psi(\uu_{f''\circ f'})}p',\alpha'\circ (\alpha''\cdot_{f'}p')\big). 
    $$
    
     \item a 2-cell  
        
\[\begin{tikzcd}
	{(A,a)} && {(B,b)}
	\arrow[""{name=0, anchor=center, inner sep=0}, "{(f,p,\alpha)}", curve={height=-12pt}, from=1-1, to=1-3]
	\arrow[""{name=1, anchor=center, inner sep=0}, "{(g,q,\beta)}"', curve={height=12pt}, from=1-1, to=1-3]
	\arrow["{(\rho,\gamma)}", shorten <=3pt, shorten >=3pt, Rightarrow, from=0, to=1]
\end{tikzcd}\]
 is a pair $(\rho,\gamma)$ consisting of a 2-cell $\rho\colon f\Rightarrow g$ of $\C$ and a map $\gamma$ in $\tP_{\psi(g)}$ which we further specify. Let us first define a map $\tP_\rho\colon \tP_{\psi(f)}\to \tP_{\psi(g)}$. We consider $\rho$ as a triangle $\bar{\rho}\colon f\circ \uu\Rightarrow g$, whose fiber is $\psi(\bar{\rho})\colon \psi(g) \to \psi(f)$ and we have $\psi(\psi(\bar{\rho}))=\psi(\uu)=v_c$ for some $c\in \pi_0(\C)$. We~define $\tP_{\rho}(p):=p\cdot_{\psi(\bar{\rho})}e_c$ and $\gamma$ is any map $\gamma\colon\tP_{\rho}(p)\to q$ in $\tP_{\psi(g)}$, such that following diagram commutes in $\tP_A$.
\[\begin{tikzcd}
	b\cdot_{g} \tP_\rho(p) & b\cdot_f p\\
	 b\cdot_{g} q& a 
	\arrow["{=}", from=1-1, to=1-2]
	\arrow["b\cdot_g \gamma"', from=1-1, to=2-1]
    \arrow["\beta", from=2-1, to=2-2]
    \arrow["\alpha", from=1-2, to=2-2]
 \end{tikzcd}\]
               \end{itemize}
             
The augmented nerve $\N\C_\tP$ consists of:
\begin{itemize}
    \item $(\N\C_\tP)_{-1}=\pi_0(\C_\tP)$, the set connected components of $\C_\tP$. This is in bijection with connected components of $\C$ via the maps $$(A,a) \xrightarrow{(u_A,a,\uu_a)} (v_{\pi(A)},e_{\pi(A)})$$ in $\C_\tP$, and also in bijection with the set of units $e_c$ of the operad $\tP$. We shall use $$(\N\C_\tP)_{-1}=\pi_0(\C).$$
    \item $(\N\C_\tP)_{0}=(\C_\tP)_0$ consists of pairs $(A\in \C_0,a\in \tP_A)$,
    interpreted as segments $$\psi(A) \xrightarrow{(A,a)}\pi(A).$$
    This gives $\phi(A,a) :=\psi(A)$.
    \item $(\N\C_\tP)_{1}=(\C_\tP)_1$ consists of maps $(f, p,\alpha )$, with $A\xrightarrow{f}B$ and \hbox{$b\cdot_f p \xrightarrow{\alpha} a$}, interpreted as triangles 
\[\begin{tikzcd}
	\psi(A) && \psi(B) \\
	& \pi(A)=\pi(B)
	\arrow[""{name=0, anchor=center, inner sep=0}, "{\big(\psi(f),p\big)}", from=1-1, to=1-3]
	\arrow["{(A,a)}"', from=1-1, to=2-2]
	\arrow["{(B,b)}", from=1-3, to=2-2]
	\arrow["{(f, p,\alpha)}"{description}, draw=none, from=0, to=2-2]
\end{tikzcd}\]
This gives $\phi(f,p,\alpha) :=\big(\psi(f),p\big)$.

\item $(\N\C_\tP)_{2}$ consists of lax triangles in $\C_\tP$, that is 2-cells $$(\rho,\gamma)\colon (g,q,\beta)\circ(f,p,\alpha)\Rightarrow (h,r,\delta)$$ which means
\begin{itemize}
\item$\rho\colon(g\circ f)\Rightarrow h$ is a 2-cell in $\C$,
\item$\gamma\colon \tP_{\rho}(q\cdot_{\psi(\uu_{g\circ f})}p)\to r$ is a map in $\tP_{\psi (h)}$, and
\item
$
\delta\circ(c\cdot_h\gamma)=\alpha\circ(\beta\cdot_f p).
$
\end{itemize}
\begin{lemma}
Every element $(\rho,\gamma)$ of $\N\C_\tP$ can be presented as a 3-simplex
$$\big((g,q,\beta),(h,r,\delta),(f,p,\alpha),(\psi(\rho),p,\gamma)). $$
\end{lemma}
\pf
One has to check that $$q\cdot_{\psi(\rho)}p= q\cdot_{\psi(\rho)}(p\cdot_{\uu_f} e_{\pi A}) = (q\cdot_{\psi(\uu_{g\circ f})}p)\cdot_{\psi(\bar{\rho})} e_{\pi A} =\tP_{\bar{\rho}}(q\cdot_{\psi(\uu_{g\circ f})}p),$$ which follows from the associativity of $\tP$ and the following  3-simplex in $\C$.
\[\begin{tikzcd}
	{A} & {A} \\
	{C} & {B}
	\arrow["{=}", from=1-1, to=1-2]
	\arrow[""{name=0, anchor=center, inner sep=0}, "{h}"', from=1-1, to=2-1]
	\arrow["{f}", from=1-2, to=2-2]
	\arrow["{g}", from=2-2, to=2-1]
	\arrow["{\uu_{g\circ f}}"{description},pos=0.55, draw=none, from=2-2, to=1]
    \arrow[""{name=1, anchor=center, inner sep=0}, "{g\circ f}"{description},  from=1-2, to=2-1]
	\arrow["{\bar{\rho}}"{description},pos=0.35, draw=none, from=1-1, to=1]
\end{tikzcd}
=
\begin{tikzcd}
	{A} & {A} \\
	{C} & {B}
	\arrow["{=}", from=1-1, to=1-2]
	\arrow[""{name=0, anchor=center, inner sep=0}, "h"', from=1-1, to=2-1]
	\arrow[""{name=1, anchor=center, inner sep=0},"f"{description}, from=1-1, to=2-2]
	\arrow["f", from=1-2, to=2-2]
	\arrow["g", from=2-2, to=2-1]
	\arrow["{\uu_f}"{description},pos=0.35, draw=none, from=1-2, to=1]
	\arrow["\rho"{description},pos=0.4, draw=none, from=2-1, to=1]
\end{tikzcd}
\]
\epf
The lemma gives $\phi (\rho,\gamma):=\big(\psi(\rho),p,\gamma\big)$.
\item
    The maps $u$ are given by
    $$u_c=(v_c,e_c),$$
    $$u_{(A,a)}=(v_A,a,\uu_{a}),$$
    $$u_{(f, p,\alpha )}=(v_f,\alpha).$$
\end{itemize}
We leave to the reader the specification of $(\NC_\tP)_3$ and the maps $\phi\colon (\NC_\tP)_3\to (\NC_\tP)_2$ and  $u\colon (\NC_\tP)_2\to (\NC_\tP)_3$, together with checking the conditions of Proposition~\ref{proposition:unary_operadic_2-category_B}.  
\begin{definition}\label{definition:G-construction2}
For any (non-1-connected) categorical $\bO$-operad we define its operadic Grothendieck construction explicitly in the very same way as for 1-connected operads above. The only difference from the construction for 1-connected operads is that there can be more objects and morphisms above the objects $v_c\in \C$, that is of form $(v_c,a)$.
     \end{definition}
    \begin{example}
    A strict monoidal category $\V$ is a categorical $\odot$-operad, for $\odot$ the terminal unary operadic (1-)category of Example~\ref{example:odot}. Its operadic Grothendieck construction $\int_\odot \V$ is isomorphic to $\V_\triangle$ of Example~\ref{example:monoidal_categories}. The assignment $\V \mapsto \int_\odot \V$ is a fully faithful functor $$\int_\odot\colon \StrMonCat\to\TwoUnOpCat.$$ In Section~\ref{section:left_adjoint_*} we describe its left adjoint.
    \end{example}
    \begin{example}
    For the operadic 1-category $\Bq{A}$ of unary bouquets on a set $A$, a~$\Bq{A}$-operad $\tP$ is a 2-category with the set of objects $A$ and $\int_{\Bq{A}} \tP \cong (\N\tP,\D\tP)$. 
    \end{example}
    \begin{example}
        If $\bO$ is a unary operadic 1-category, and $\tP$ a categorical $\bO$-operad which is discrete in each component (i.e.~an operad in $\Set$), we recover the discrete operadic Grothendieck construction of \cite{Batanin_Markl:duoidal}.     
        \end{example}

\section{Operadic fibrations}
Operadic fibrations (with splitting) are operadic functors lying in the essential image of the functor~$\int$ of the preceding section. We characterize them by certain lifting properties.  A prototypical operadic fibration is thus the projection $\pi_\tP\colon \int_\bO \tP\to \bO$ from the previous section.
\begin{definition}
Let $\bO=(X,\C)$ be a unary operadic 2-category. The objects $u_c \in \mathscr{U}$, i.e.~objects lying in the image of $u\colon \pi_0\C \to \C_0$, are called trivial.
A morphism $g\colon y\to x$ of a unary operadic 2-category is a quasibijection if for any map $f\colon z\to y$, $$\phi(g\circ f\xRightarrow{\uu_{g\circ f}}(g \circ f))=u_{\phi(f)}.$$
\end{definition}

It follows that $g$ has a trivial fiber $u_{\phi(x)}=u_{\phi(y)}$, which recovers 1-categorical quasibijections of \cite[Section~1.1]{Batanin_Markl:kodu2021}.
\begin{proposition}\label{proposition:quasibijections} Let $\bO=(X,\C)$ be a unary operadic 2-category and $\tP$ a categorical $\bO$-operad.
A~morphism $G=(g,q,\beta)\colon (B,b)\to (C,c)$ of $\int_\bO \tP$ is a quasibijection if and only if $g$ is a quasibijection in $\bO$ and $q=e_{\phi(B)}$, the unit of the operad $\tP$. 
\end{proposition}
\pf
Let $(f,p,\alpha)\colon (A,a)\to (B,b)$ be a morphism of $\int_\bO \tP$ and $g\colon B\to C$ a quasibijection in $\bO$. We check that $$\phi\big(\uu_{(g,e_{\phi(B)},\beta)\circ (f,p,\alpha)}\big) = \phi(\uu_{g \circ f},\uu_p) =\big(\phi(\uu_{g \circ f}),p,\uu_p\big)= (u_{\phi(f)},p,\uu_p) = u_{\phi(f,p)}.$$
Conversely, let $(g,q,\beta)$ be a quasibijection. We get the equation 
$$(\phi(\uu_{(g\circ f)}),q\cdot p,\uu_{q\cdot p})= \phi(\uu_{g\circ f},\uu_{q\cdot p})=\phi(\uu_{(g,q,\beta)\circ (f,p,\alpha)})=u_{\phi(f,p)}=(u_{\phi(f)},p,\uu_{p})$$
which holds for any $f\colon A\to B$ and $p\in \tP_{\phi(f)}$. From the first components we see immediately that $g$ is a quasibijection in $\bO$. By evaluating $f=\uu_B$ and $p=e_{\phi(B)}$ we get $q=e_{\phi(B)}$.   
\epf
\begin{definition}
    For $k\geq 0$, $\Delta^k$ denotes the simplicial set $\Set(-,[k])$. A $(3,2)$-horn is the sub-simplicial set $\iota\colon\Lambda_2^3 \hookrightarrow \Delta^3$, which is the union of all faces except the face $(013)$. A~$(3,2)$-horn $\lambda=(d_0\lambda,d_1\lambda,d_3\lambda)$ in a simplicial set $X$ is a morphism $\lambda\colon \Lambda_2^3 \to X$. A filler of a (3,2)-horn $\lambda\colon \Lambda_2^3 \to X$ is a 3-simplex $\sigma \colon \Delta^3 \to X$, whose restriction along the inclusion $\iota$ is $\lambda$.
\end{definition}
\begin{definition}
    Let $p\colon \bP \to \bO$, be a morphism of operadic 2-categories with $$\bP=(X,\C), \bO=(Y,\D),$$ and let $f\in X_2$. We say that $f$ is operadic $p$-cartesian if for any (3,2)-horn $\lambda$ in $X$ with $d_0\lambda= f$ and any filler  $\sigma$ of $p\lambda$ in $Y_2$, there exists a unique filler $\Tilde{\sigma}$ of $\lambda$, such that $p\Tilde{\sigma}=\sigma$.
\end{definition}
\begin{definition}\label{definition:fibration}
   A morphism $p\colon \bP \to \bO$ is an operadic fibration if it induces an epimorphism  $\pi_0\bP \twoheadrightarrow \pi_0\bO$ on the sets connected components, and for any $a,b \in X_1$ and $g\in Y_2$ with $d_0g=p(a)$ and $d_2g=p(b)$, there exists an operadic $p$-cartesian lift $\Tilde{g}$ of $g$ in~$\C$. An operadic fibration is split if there is a fixed choice of $\Tilde{g}$ for every $a,b,g$, such that the following two conditions hold. The~choice of $\Tilde{g}$ will be denoted by $\ell(a,b,g)$.
   \begin{itemize}
       \item For every $x \in X_1$, $\ell(u_{\pi(x)},x,u_{p(x)}) = u_x$ and $\ell(x,u_{\phi (x)},\uu_{p(x)}) = \uu_x$. 
       \item For every 3-simplex $\sigma_{0123}\in Y_3$  with $f_{01}=p(x)$, $f_{12}=p(y)$, and $f_{23}=p(z)$ (cf.~Section~\ref{section:Nerve_and_dec}), for~some $x,y,z \in X_1$, 
       $$\ell(d_1\ell(z,y,\alpha_{123}),x,\alpha_{013})=\ell(z,d_1 \ell(y,x,\alpha_{012}),\alpha_{023}).$$
   \end{itemize}
\end{definition}
\begin{proposition}\label{proposition:iso_on_components}
A split operadic fibration $p\colon \bP \to \bO$ induces an isomorphism on connected components.
\end{proposition}
\pf Let $c,d\in \pi_0(\D)$ with $p(c)=p(d)$. Then $p(u_d)=p(u_c)=u_{p(c)}$. By Proposition~\ref{proposition:unary_operadic_2-category_B}, \eqref{item:15}, $u_{\phi(u_c)}=u_c$. The splitting conditions give $\ell(u_d,u_c,\uu_{u_{p(c)}})=\uu_{u_c}$. Comparing fibers of these two maps we see $u_c=u_d$, and hence $c=d$.   \epf
\begin{corollary}\label{corollary:unique_trivial}
Let $p\colon \bP \to \bO$ be a split operadic fibration. For any trivial object $v$ of $\bO$ there is a unique trivial object $u$ in $\bP$ above $v$.
\end{corollary}
\begin{proposition}
For a unary operadic 2-category $\bO$ and a categorical $\bO$-operad $\tP$, the operadic functor $\pi_\tP \colon \int_\bO \tP\to \bO$ is a split operadic fibration.
\end{proposition}
\proof
For objects $(A,a)$ and $(B,b)$ of $\int_\bO \tP$, and morphism $f\colon A\to B$ in $\bO$, the operadic $\pi_\tP$-cartesian lift of $f$ is $(f,b,\uu_{a\cdot_f b})$. It is straightforward to check the conditions. 
\endproof
\begin{theorem}\label{theorem:fibrations}
 There is a one-to-one correspondence of categorical $\bO$-operads and split operadic fibrations over a unary operadic 2-category $\bO$.  
\end{theorem}

\pf For a split operadic fibration $p\colon \bP\to\bO$ we define a categorical $\bO$-operad~$\tP$. Let $\bP=(Y,\D) $ and $\bO=(X,\C)$. For $x\in \C$, let $\tP_x:=p_{\mathrm{qb}}^{-1}(\uu_x)$ be the subcategory of $\D$ of all objects which are mapped to $x$ and quasibijections which are mapped to $\uu_x$. For a morphism $g$ of $\C$, that is $g\in X_2$, the operad multiplication $\cdot_g\colon \tP_{d_0g}\x \tP_{d_2g}\to \tP_{d_1g}$ is defined using the $p$-cartesian lifts: for $a\in \tP_{d_0g}$ and $b\in \tP_{d_2g}$, $a\cdot_g b:= d_1 \ell(a,b,g)$. Let $\alpha\colon a' \to a''$ in $\tP_{d_0g}$ and $\beta\colon b' \to b''$ in $\tP_{d_2g}$. Since $\alpha$ is a quasibijection, $\phi(\uu_{\alpha\circ \ell(a',b',g)})=u_{b'}$. There is a $(3,2)$-horn $\lambda = (\ell(a'',b'',g),\alpha\circ \ell(a',b',g),\beta)$, and $p\lambda$ has a filler $\sigma=(g,g,\uu,\uu)$ in~$\C$. Since $\ell(a'',b'',g)$ is $p$-cartesian, there is a unique filler $\Tilde{\sigma}$ of $\lambda$ in $\D$. This defines $\alpha\cdot_g \beta:=d_2\Tilde{\sigma}$. 
 For $c\in \pi_0(\C)$, by Corollary~\ref{corollary:unique_trivial}, there is a unique trivial object $e_c \in \tP_{u_c}$, which provides the unit of $\tP$.
 The functoriality, unitality and associativity of the operad $\tP$ is ensured by the splitting conditions of Definition~\ref{definition:fibration} and the uniqueness of horn fillers. Using Proposition~\ref{proposition:quasibijections}, it is straightforward to check that this construction is inverse to~$\int_{\bO}$.
 \endproof

\begin{example}
     Let $\bO=\odot$, the terminal unary operadic 1-category. Theorem~\ref{theorem:fibrations} provides a (possibly new) characterization of strict monoidal categories. Namely, a strict monoidal category is a unary operadic 2-category $\bM=(X_{\bM},\C_{\bM})$ with the following lifting property. For any two objects $a,b\in \C_{\bM}$, there exists a universal object $a\cdot b$ and a universal morphisms $\ell(a,b)\colon a\cdot b \to a$ in $\C_{\bM}$ with $\phi (\ell(a,b))=b$, and such that for any pair of maps $k\colon c\to a$ and $h\colon \phi (k) \to b$, there is a unique morphism  $g\colon c \to  a\cdot b$ with $\phi(g) = \phi(h)$ and $(f,g,k,h)\in (\N\C_{\bM})_3$. The splitting provides a unique trivial object $e$ of $\bM$ and conditions $\ell(e,x)=u_x$, $\ell(x,e)=\uu_x$, and $$\ell(d_1\ell(a,b),c)=\ell(a,d_1\ell(b,c)).$$ By looking at domains of these maps we get $e\cdot x=x= x\cdot e$ and $(a\cdot b)\cdot c= a\cdot(b\cdot c)$.

 \end{example}

   \section{Left adjoint to $\int_\odot$}\label{section:left_adjoint_*}
 Recall the terminal unary operadic category $\odot$ of Example~\ref{example:odot}, whose categorical operads are strict monoidal categories. Considering $\odot$ as a 2-category with the only identity 2-cell, it is a terminal operadic 2-category.
In this section we construct a left adjoint to the functor $\int_\odot\colon \StrMonCat \to \TwoUnOpCat$. First we recall a few familiar adjunctions.\\


Let $\Delta_{\leq n}$ denote the full subcategory of $\Delta$, which consists only of objects $[0], \ldots, [n]$. An~$n$-truncated simplicial set is a presheaf $\Delta_{\leq n}^{\op}\to \Set$. The category of $n$-truncated simplicial sets is denoted by $\sSetn$. The inclusion $\Delta_{\leq n} \subseteq\Delta$  induces an adjunction
    \[\begin{tikzcd}
            \sSetn \arrow[r, shift left=1ex, "\textrm{cosk}_n"{name=U}] & \sSet\arrow[l, shift left=.5ex, "{\textrm{tr}_n}"{name=F}]
            \arrow[phantom, from=F, to=U, , "\scriptscriptstyle\boldsymbol{\top}"],
        \end{tikzcd}\]
    where $\textrm{tr}_n$ is the restriction along the inclusion, and $\textrm{cosk}_n$ is given by right Kan extension.
    
   For $S\in \Set$, $\coll{S}{\Set}$ denotes the category of functors $$\coprod_{n\geq 0} \big(S^{\x n}\x S\big) \xrightarrow{Q} \Set,$$ and $\textrm{Multicat}_S$ denotes the category of multicategories with the set of objetcs $S$. 
   An $S$-collection~$Q$ is just a collection of sets $Q(a_1\cdots a_k;a_0), a_i \in S, k\geq 0.$
   There is the free-forgetful adjunction   
    \[\begin{tikzcd}
            \Multi_S \arrow[r, shift left=1ex, "U"{name=U}] & \coll{S}{\Set}\arrow[l, shift left=.5ex, "F"{name=F}]
            \arrow[phantom, from=F, to=U, , "\scriptscriptstyle\boldsymbol{\top}"].
        \end{tikzcd}\]
        The multimorphisms of $FQ$ are planar rooted trees with inner vertices labeled by elements of the collection $Q$.
        
        Last, we recall the adjunction 
        \[\begin{tikzcd}
            \StrMonCat \arrow[r, shift left=1ex, "R"{name=U}] & \Multi \arrow[l, shift left=.5ex, "L"{name=F}]
            \arrow[phantom, from=F, to=U, , "\scriptscriptstyle\boldsymbol{\top}"].
        \end{tikzcd}\]
For~$\M\in \StrMonCat$, the multicategory $R\M$ has $ R\M(a_1\cdots a_k;a_0)=\M(a_1\ox\cdots\ox a_k,a_0)$. For $\tM\in \Multi$, the objects of $L\tM$ are lists of objects of $\tM$. A morphism of lists $$(a_1,\cdots,a_k) \to (b_1,\cdots,b_l)$$ is an $l$-tuple tuple of multimorphisms of $\tM$. The monoidal product is concatenation of lists and the monoidal unit is the empty list.

We first establish an adjunction
  \[\begin{tikzcd}
            \StrMonCat \arrow[r, shift left=1ex, "\Psi"{name=U}] & \sSett \arrow[l, shift left=.5ex, "\Phi"{name=F}]
            \arrow[phantom, from=F, to=U, , "\scriptscriptstyle\boldsymbol{\top}"].
        \end{tikzcd}\]

        For $\M\in \StrMonCat$ we define $\Psi \M=\textrm{tr}_3 \N\B\M$, the 3-truncated nerve of the one-object 2-category $\B\M$. 
\begin{lemma}
    Let $\C$ be a small category, $\sP\C$ the category of presheaves on $\C$, and \hbox{$G\colon \D \to \sP\C$} any functor from a cocomplete category. Let $F_0$ be a functor $F_0\colon \C \to \D$, such that for any object $s \in \C$ and $a\in \D$, there is a natural isomorphism $$\D(F_0b,a)\cong \sP\C(yb,Ga),$$ where $y\colon \C \to \sP\C$ is the Yoneda embedding. Then $F=Lan_yF_0$, the left Kan extension of $F_0$ along $y$, is left adjoint to $G$.
    \end{lemma}
\proof
Any $X\in \sP\C$ is isomorphic to the colimit of representables $X \cong \colim yb$, where the colimit is taken over the category of elements of $X$. 
There are natural isomorphisms $$\D(FX,a)\cong \D(F\colim yb,a)\cong \D(\colim Fyb,a)\cong \D(\colim F_0b,a)\cong \lim \D(F_0b,a)\cong 
$$
$$\cong \lim \sP\C(yb,Ga)\cong \sP\C(\colim yb,Ga)\cong \sP\C(X,Ga).\qedhere $$

        By the preceding lemma it suffices to define a functor $\Phi_0 \colon \Delta_{\leq3}\to \StrMonCat$ with $$\StrMonCat(\Phi_0[k],M) \cong \sSett(\Delta^k,\Psi M),$$ for $k=0,1,2,3,$ where $\Delta^k:=y[k]=\Set(-,[k])$. By the Yoneda lemma, $$\sSett(\Delta^k,\Psi M)\cong (\Psi M)_k.$$

\begin{itemize}
    \item For $k=0$, we define $\Phi_0[0]=1$, the terminal category, which is also the initial monoidal category. We see that 
    $$1\cong (\Psi M)_0\cong \StrMonCat(1,M).$$
   \item For $k=1$, we define $\Phi_0[1]=(\mathbb{N},+,0)$, the free strict monoidal category on one object. A~monoidal functor $H\colon \mathbb{N}\to M$ is determined by  an object $H(1)\in M$, hence $$M_0=(\Psi M)_1\cong \StrMonCat(\mathbb{N},M).$$
   \item For $k=2$, we consider the simplicial set $\Delta^2$, whose non-degenerate simplicies are: 
   \[\begin{tikzcd}
	x_0 && x_1. \\
	& x_2
	\arrow["f_{01}", from=1-1, to=1-3]
	\arrow[""{name=0, anchor=center, inner sep=0}, "f_{02}"', from=1-1, to=2-2]
	\arrow["f_{12}", from=1-3, to=2-2]
	\arrow["\alpha_{012}"{description},pos=0.4, draw=none, from=0, to=1-3]
\end{tikzcd}\]
   Let $S=\{f_{01},f_{12},f_{02}\}$. We define an $S$-collection $Q$ by $Q(f_{12}\,f_{01};f_{02})=\{\alpha_{012}\}$ and empty otherwise. We put \hbox{$\Phi_0[2]=LFQ.$}  Objects of the monoidal category $\Phi_0[2]$ are lists $(a_1,\ldots,a_k)$, $a_i\in \{f_{01},f_{12},f_{02}\}$, $k\geq 0$. Morphisms $(a_1,\ldots,a_k)\to (b_1,\ldots,b_l)$ are tuples $(g_1,\ldots,g_l)$, where each $g_i$ is either~$\alpha_{012}$ or the identity. In other words, $\Phi_0[2]$ is the universal strict monoidal category containing a lax triangle.
   A monoidal functor $H\colon \Phi_0[2] \to M$ is determined by three objects $H(f_{01}),H(f_{12}),H(f_{02})$ of~$M$ and a morphism $H(f_{12})\ox H(f_{01})\xrightarrow{H(\alpha_{012})} H(f_{02})$, hence
   $$(\Psi M)_2\cong \StrMonCat(\Phi_0[2],M).$$
\item For $k=3$, we consider a filled 3-simplex $\sigma_{0123}$ with boundary  
\[\begin{tikzcd}
	{x_0} & {x_1} \\
	{x_3} & {x_2}
	\arrow["{f_{01}}", from=1-1, to=1-2]
	\arrow[""{name=0, anchor=center, inner sep=0}, "{f_{03}}"', from=1-1, to=2-1]
	\arrow[""{name=1, anchor=center, inner sep=0},  from=1-2, to=2-1]
	\arrow["{f_{12}}", from=1-2, to=2-2]
	\arrow["{f_{23}}", from=2-2, to=2-1]
	\arrow["\alpha_{013}"{description},pos=0.35, draw=none, from=1-1, to=1]
	\arrow["\alpha_{123}"{description},pos=0.4, draw=none, from=2-2, to=1]
\end{tikzcd}
 \hspace{1cm}
\begin{tikzcd}
	{x_0} & {x_1} \\
	{x_3} & {x_2}
	\arrow["{f_{01}}", from=1-1, to=1-2]
	\arrow[""{name=0, anchor=center, inner sep=0}, "{f_{03}}"', from=1-1, to=2-1]
	\arrow[""{name=1, anchor=center, inner sep=0}, from=1-1, to=2-2]
	\arrow["{f_{12.}}", from=1-2, to=2-2]
	\arrow["{f_{23}}", from=2-2, to=2-1]
	\arrow["\alpha_{012}"{description},pos=0.35, draw=none, from=1-2, to=1]
	\arrow["\alpha_{023}"{description},pos=0.4, draw=none, from=2-1, to=1]
\end{tikzcd}
\]
Let $B$ be the set of its 1-simplicies, 
$$B=\{f_{ij}\,|\, 0\leq i<j\leq3\}. $$
We define a $B$-collection $R$ by 
\begin{itemize}
    \item $R(f_{jk}\,f_{ij};f_{ijk})=\{\alpha_{ijk}\,|\, 0\leq i<j<k\leq3\}$.
    \item $R(f_{23}\,f_{12}
    ,f_{01};f_{03})=\{\sigma_{0123}\}$, and empty otherwise.
\end{itemize}
Consider \hbox{$FR/\sim,$} the quotient of the free multicategory, where $\sim$ is generated by: 
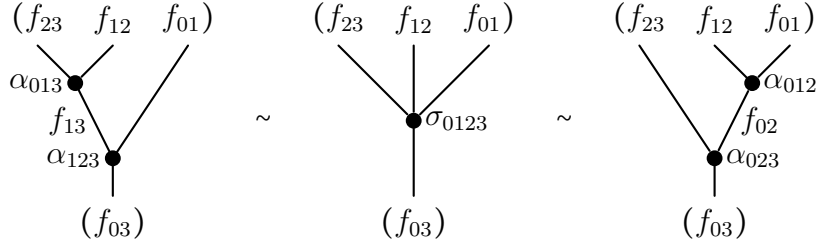
\[\begin{tikzpicture}
[thick, line cap=round,line join=round,
    x=1.0cm,y=1.0cm
    ] 
\draw (0,-0.5) node[below]{$(f_{03})$}--(0,0)node[circle,fill=black,inner sep=0pt,minimum size=6pt](p){};
 \draw  (-1,1.5)node[above]{$( f_{23}$} -- (-0.5,1)node[circle,fill=black,inner sep=0pt,minimum size=6pt](q){}; 
 \draw  (0,1.5)node[above]{$f_{12}$} -- (q) ;
 \draw  (1,1.5)node[above]{$f_{01})$} -- (0,0) ;
 \draw (-0.5,1) -- (0,0);
\draw (q) node[left] {$\alpha_{013}$};
\draw (p) node[left]{$\alpha_{123}$}; 
\draw (-0.6,0.5)node{$f_{13}$};

\begin{scope}[shift={(4,0)}]
[thick, line cap=round,line join=round,
    x=1.0cm,y=1.0cm
    ] 
\draw (0,-0.5) node[below]{$(f_{03})$}--(0,0.5)node[circle,fill=black,inner sep=0pt,minimum size=6pt](p){};
 \draw  (-1,1.5)node[above]{$(f_{23}$} -- (p);
 \draw  (0,1.5)node[above]{$f_{12}$} -- (p) ;
 \draw  (1,1.5)node[above]{$f_{01})$} -- (p) ;
\draw (p) node[right]{$\sigma_{0123}$}; 
\end{scope}

\begin{scope}[shift={(8,0)}]
[thick, line cap=round,line join=round,
    x=1.0cm,y=1.0cm
    ] 
\draw (0,-0.5) node[below]{$(f_{03})$}--(0,0)node[circle,fill=black,inner sep=0pt,minimum size=6pt](p){};
 \draw  (-1,1.5)node[above]{$(f_{23}$} --(p);
 \draw  (0,1.5)node[above]{$f_{12}$} --(0.5,1) node[circle,fill=black,inner sep=0pt,minimum size=6pt](q){};  ;
 \draw  (1,1.5)node[above]{$f_{01})$} --(q) ;
 \draw (q) --(p);
\draw (q) node[right] {$\alpha_{012}$};
\draw (p) node[right]{$\alpha_{023}$}; 
\draw (0.6,0.5)node{$f_{02}$};
\end{scope}
\draw (2,0.5)node{$\sim$};
\draw (6,0.5)node{$\sim$};
\end{tikzpicture}\]

We put $\Phi[3]=L(FR/\sim)$. Equivalently, the monoidal category $\Phi[3]$ consists of objects $(a_1,\ldots,a_k)$, $a_i\in B$, $k\geq 0$, and morphisms $(a_1,\ldots,a_k)\to (b_1,\ldots,b_l)$ are tuples $(g_1,\ldots,g_l)$, where each $g_i$ is either an identity on $f_{ij}\in B$, a non-degenerate 2-simplex $\alpha_{ijk}$, or the 3-simplex~$\sigma_{0123}$. An example of composition of morphisms is:
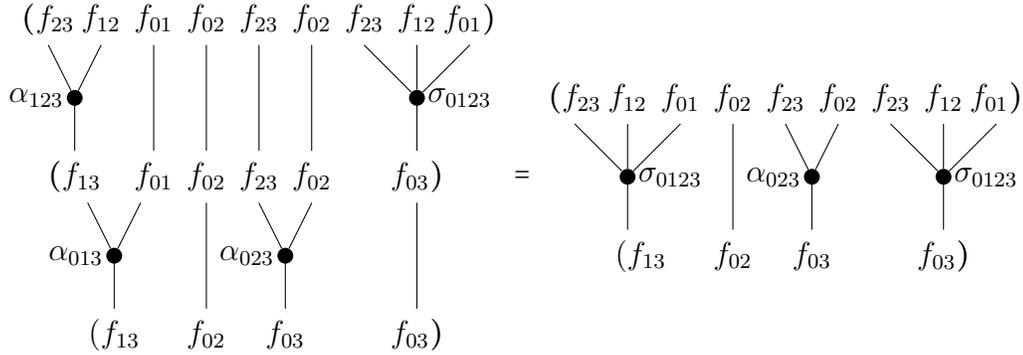
\[\begin{tikzpicture}
[line cap=round,line join=round,
    x=0.7cm,y=0.7cm
    ]
    \begin{scope}[shift={(0,3.5)}]
    \draw (1,0)node{$(f_{23}$};
    \draw (2,0)node{$f_{12}$};
    \draw (3,0)node{$f_{01}$};
    \draw (4,0)node{$f_{02}$};
    \draw (5,0)node{$f_{23}$};
    \draw (6,0)node{$f_{02}$};
    \draw (7,0)node{$f_{23}$};
    \draw (8,0)node{$f_{12}$};
    \draw (9,0)node{$f_{01})$};
    \end{scope}
    \begin{scope}[shift={(0,1)}]
    \tikzalpha(1.5;123)
    \tikzf(3)
    \tikzf(4)
    \tikzf(5)
    \tikzf(6)
    \tikzsigma(8)
    \end{scope}
    \begin{scope}[shift={(0,0.5)}]
    \draw (1.5,0)node{$(f_{13}$};
    \draw (3,0)node{$f_{01}$};
    \draw (4,0)node{$f_{02}$};
    \draw (5,0)node{$f_{23}$};
    \draw (6,0)node{$f_{02}$};
    \draw (8,0)node{$f_{03})$};
    \end{scope}
    \begin{scope}[shift={(0,-2)}]
    \tikzalpha(2.25;013)
    \tikzf(4)
    \tikzalpha(5.5;023)
    \tikzf(8)
    \end{scope}
    \begin{scope}[shift={(0,-2.5)}]
    \draw (2.25,0)node{($f_{13}$};
    \draw (4,0)node{$f_{02}$};
    \draw (5.5,0)node{$f_{03}$};
    \draw (8,0)node{$f_{03})$};
    \end{scope}
    \draw (10,0.5)node{$=$};
  \begin{scope}[shift={(10,2)}]
   \draw (1,0)node{$(f_{23}$};
    \draw (2,0)node{$f_{12}$};
    \draw (3,0)node{$f_{01}$};
    \draw (4,0)node{$f_{02}$};
    \draw (5,0)node{$f_{23}$};
    \draw (6,0)node{$f_{02}$};
    \draw (7,0)node{$f_{23}$};
    \draw (8,0)node{$f_{12}$};
    \draw (9,0)node{$f_{01})$};
    \end{scope}
     \begin{scope}[shift={(10,-0.5)}]
    \tikzsigma(2)
    \tikzf(4)
    \tikzalpha(5.5;023)
    \tikzsigma(8)
    \end{scope}
    \begin{scope}[shift={(10,-1)}]
    \draw (2.25,0)node{($f_{13}$};
    \draw (4,0)node{$f_{02}$};
    \draw (5.5,0)node{$f_{03}$};
    \draw (8,0)node{$f_{03})$};
    \end{scope}
\end{tikzpicture}\]
   A monoidal functor $H\colon \Phi[3] \to M$ is determined by six objects $H(f_{ij})\in M$, four morphisms $$H(f_{jk})\ox H(f_{ij})\xrightarrow{H(\alpha_{ijk})} H(f_{ik}),$$ and a morphism $$H(f_{23})\ox H(f_{12})\ox H(f_{01})\xrightarrow{H(\sigma_{0123})} H(f_{03}),$$ which satisfy 
   $$H(\alpha_{123})\circ \big(H(\alpha_{013})\ox \uu_{H(f_{01})}\big)=H(\sigma_{0123})= H(\alpha_{023})\circ\big(\uu_{H(f_{23})}\ox H(\alpha_{012})\big).$$
   This gives $$(\Psi M)_3\cong \StrMonCat(\Phi[3],M).$$
  \end{itemize}
By a $k$-simplex in a monoidal category $\M$ we mean a $k$-simplex in the 2-category $\B\M$ (cf.~Definition~\ref{definition:duskin_nerve}). Denote by $\sigma_l$ the generating $l$-simplex of $\Phi_0[l]$.
Let $k,l\leq 3$ and $f\colon [k]\to [l]$.
By~the~above $$\StrMonCat(\Phi_0[k],\Phi_0[l])\cong \{\text{$k$-simplices of $\Phi_0[l]$}\},$$ i.e.~any monoidal functor $\Phi_0[k]\to\Phi_0[l]$ is determined by a $k$-simplex in $\Phi_0[l]$. In particular, $\Phi_0 f$ is determined by the $k$-simplex $\sigma_l\circ f$.
This completes the functor $$\Phi_0\colon \Delta_{\leq3}\to \StrMonCat,$$ such that the above isomorphisms are natural, which finishes the construction of the adjunction
 \[\begin{tikzcd}
            \StrMonCat \arrow[r, shift left=1ex, "\Psi"{name=U}] & \sSett .\arrow[l, shift left=.5ex, "\Phi"{name=F}]
            \arrow[phantom, from=F, to=U, , "\scriptscriptstyle\boldsymbol{\top}"]
        \end{tikzcd}\]
        
        Since $\N\B\M$ is 3-coskeletal, $\N\B\M\cong \textrm{cosk}_3( \textrm{tr}_3  (\N\B\M))$, and from the adjunction $$\textrm{tr}_3\dashv\textrm{cosk}_3$$ we have
        $$\sSet(X,\N\B\M)\cong \sSett(\textrm{tr}_3X,\textrm{tr}_3\N\B\M)\cong \StrMonCat(\Phi \textrm{tr}_3X, M).$$
        Thus, for a unary operadic 2-category $\bO=(X,\C)$, we have the adjunction $\Phi\mathrm{tr}_3 \dashv \int_\odot$,
        $$\TwoUnOpCat\big(\bO,\int_\odot\M)\big):=\sSet(X,\N\B\M)\cong \StrMonCat(\Phi\mathrm{tr}_3 X, \M).$$    

          \begin{example} Let $\C$ be a 1-category (viewed as a 2-category with only trivial \hbox{2-cells}) and \hbox{$\bC=(\N\C,\D\C)$} its associated unary operadic 2-category of Example~\ref{example:monoidal_categories}. The~monoidal category \hbox{$\M_\C:=\Phi\mathrm{tr}_3\N\C$} has objects lists $(f_1,\ldots,f_k)$, $f_i \in \C_1,k\geq 0$, and morphisms are generated~by $$(f_1,\ldots, g,h,\ldots, f_k)\xrightarrow{\textrm{compose}} (f_1,\ldots, g\circ h,\ldots, f_k).$$
          The monoidal product is concatenation and the monoidal unit is the empty list.
        \end{example}

        \section{Closing remarks and further directions}
        
There is a straightforward generalization of our theory to unary operadic bicategories and categorical pseudo-operads (operads with relaxed associativity and unitality conditions). We expect that these correspond to non-splitting fibrations. Working in this generalized framework, one could work directly with non-strict monoidal categories. We also hope to resolve the 1-connectedness restriction. The proofs might become more complicated however. 

A direct continuation of the work is then to establish the 2-categorical (and also bicategorical) framework for operadic categories of general arities. This is a work in progress. Unfortunately, there is no such clear characterization using simplicial sets as in the unary case, which complicates the situation. We believe that a good path is to follow Proposition~\ref{proposition:unary_operadic_2-category_C}. The proposed definition~is:
\begin{definition}
An operadic bicategory is a bicategory $\bO$ equipped with 
\begin{itemize}
\item a cardinality 2-functor $\card{-}\colon\bO\to \mathscr{S}$ (necessarily strict) to the skeletal category of finite sets,  
\item lax functors $\phi_x\colon \bO/x\to\bO^{\card{x}}$, for every object $x$ of $\bO$, and 
\item chosen local lali-terminal objects $u_c$, in each connected component $c$ of $\bO$ whose cardinality is not constantly zero, 
\item such that analogues of the axioms of operadic 1-categories hold.
\end{itemize}
\end{definition}
 When the cardinality functor is constantly zero, the theory should recover ordinary categories. In this case the operads are presheaves and operadic fibrations are categorical fibrations. 

Let $F\colon \bO \to \bP$ be an operadic functor. It induces a restriction functor on the categories of operads $$F^*\colon \oper{\bP}{\Cat}\to\oper{\bO}{\Cat}.$$ A non-trivial problem is to find a formula for its left adjoint. Evaluating this on the terminal $\bP$-operad we possibly get a left adjoint to the operadic Grothendieck construction for general $\bO$.

\section*{Acknowledgements}
This work was supported by RVO: 67985840 and Praemium Academi\ae of M.~Markl; and by MUNI/A/1569/2024: Specific research - support for student projects.
\bibliographystyle{plain}

\end{document}